\documentclass{article}
\usepackage{amssymb, verbatim, multirow, xcolor, graphicx, caption, subcaption, sidecap}
\usepackage{hyperref}
\usepackage{amsmath}
\usepackage[utf8]{inputenc}
\usepackage[english]{babel}

\newcounter{defin}[section]

\newcounter{remark}

\def\XXint#1#2#3{{\setbox0=\hbox{$#1{#2#3}{\int}$ }
		\vcenter{\hbox{$#2#3$ }}\kern-.6\wd0}}

\newcommand{\qed}{\hfill$\Box$}
\newtheorem{theorem}{Theorem}[section]
\newtheorem{proposition}[theorem]{Proposition}
\newtheorem{lemma}[theorem]{Lemma}

\newtheorem{corollary}[theorem]{Corollary}

\title{Pseudoconvex submanifolds in K\"ahler 4-manifolds}
\author{Brian Weber}
\date{April 26, 2021}

\begin{document}

\title{Pseudoconvex submanifolds in K\"ahler 4-manifolds}
\author{Brian Weber}
\date{August 1, 2022}

\maketitle

\begin{abstract}
On K\"ahler 4-manifolds, not necessarily compact or of finite topological type, we obtain relationships between the fundamental group of compact embedded Levi-flat or pseudoconvex submanifold and the fundamental group of the ambient manifold $M^4$.
When a Levi-flat submanifold $V^3$ has finite fundamental group then $\pi_1(M^4)=\iota_*\pi_1(V^3)$; when a non-separating pseudoconvex submanifold $V^3$ has finite fundamental group, then $\pi_1(M^4)=\iota_*\pi_1(V^3)\rtimes\mathbb{Z}$.
As applications, if a K\"ahler manifold (compact or not) has an embedded holomorphic $\mathbb{P}^1$ of positive self-intersection, it must intersect all other holomorphic $P^1$ of non-negative self-intersection, the fundamental group of $M^4$ is trivial, and no ALE or ALF ends exist.
If a Levi-flat submanifold and an embedded holomorphic $\mathbb{P}^1$ of positive self-intersection both exist, they intersect.
The total number of ALE plus ALF ends is zero or one regardless of what other kinds of ends exist.
We provide examples, such as a 2-ended scalar-flat K\"ahler metric conformal to the Taub-NUT.
\end{abstract}

\section{Introduction}

We show that pseudoconvex or Levi-flat submanifolds of finite fundamental group within a K\"ahler 4-manifold (not necessarily closed or of finite topological type) force strong restrictions on the fundamental group on the ambient manifold.
As applications, we show the existence of any such a submanifold forces restrictions on the geometry of its ends, on whether any other pseudoconvex submanifolds can exist, on whether embedded holomorphic $\mathbb{P}^1$ submanifolds of positive self-intersection can exist, and restrictions on the fundamental group of any K\"ahler 4-manifold that contains such a $\mathbb{P}^1$.

We consider K\"ahler manifolds-with-boundary $(M^4,J,g)$ that satisfy the following pseudconvexity condition ($P$):
\begin{itemize}
	\item[($P$)]\hfill\parbox{.98\linewidth}{
		$M^4$ is a topologically closed manifold or manifold-with-boundary, not necesarily compact.
		If $M^4$ has boundary, every boundary component is of class $C^2$, pseudoconvex with respect to the outward pointing normal, and compact.
		There may be finite or infinitely many boundary components.
	}
\end{itemize}
Sometimes we will say ``pseudoconvex K\"ahler manifold'' to mean a manifold that satisfies ({\it{P}}).
Notably $M^4$ need not be compact or even have finite topological type.
Embedded pseudoconvex submanifolds $\iota:V^3\hookrightarrow{}M^4$, except where otherwise stated, will be assumed to have the following finiteness condition ($F$):
\begin{itemize}
	\item[($F$)]\hfill\parbox{.98\linewidth}{
		The submanifold $\iota:V^3\hookrightarrow{}M^4$ is embedded, compact, of class at least $C^2$, and is Levi-flat or pseudoconvex with respect to some choice of normal.
		Either $V^3$ intersects no point of the boundary, or it is a boundary component.
		If it is a boundary component, we choose the outward pointing normal.
		Finally $\pi_1(V^3)$ is a finite subgroup.
	}
\end{itemize}
Theorem \ref{ThmLeviFlats} shows that when $M^4$ satisfies ($P$) and $V^3\subset{}M^4$ is a Levi-flat submanifold satisfying ($F$) then in fact $\iota_*\pi_1(V^3)=\pi_1(M^4)$.
Further, any such submanifold $V^3$ is separating, and if there are two Levi-flat submanifolds $V_0^3$, $V_1^3$ of class $C^2$, and either one of them has finite fundamental group, they intersect.
If $V^3\subset\partial{}M^4$ is a boundary component it is the only boundary component.
If $V^3\subset{}M^4$ is in the interior and has finite fundamental group, then $M^4$ has no boundary.

Theorem \ref{ThmPseudoconvexSubs} weakens the Levi-flat condition on the submanifold $V^3$ to the requirement that $V^3$ be pseudoconvex rather than Levi-flat.
If $M^4$ is compact then $V^3$ is separating, and the component on which $V^3$ forms a pseudoconvex boundary has finite $\pi_1(M^4)$.
If $V^3$ is not separating then necessarily the manifold is non-compact, and its fundamental group is $\pi_1(M^4)=\iota_*(\pi_1(V^3))\rtimes\mathbb{Z}$.

Using the techniques developed to prove Theorems \ref{ThmLeviFlats} and \ref{ThmPseudoconvexSubs}, we prove a related result, Proposition \ref{PropEssSepPseudo}, dealing with the case that two pseudoconvex submanifolds $V_0^3$ and $V_1^3$ exist within $M^4$.
We cannot retain the very strong conclusion that they intersect---consider the case of small balls around distinct points in $\mathbb{C}^2$.
But we prove that no component of $M^4\setminus\{V_0^3\cup{}V_1^3\}$ cannot have pseudoconvex copies of \textit{both} $V_0^3$ and $V_1^3$ on its boundary (whether or not it might also have pseudo-concave copies of $V_0^3$ or $V_1^3$ on its boundary).
This proposition has some strong consequences.
For example a pseudoconvex K\"ahler manifold (compact or not) can have at most one embedded holomorphic $\mathbb{P}^1$ of positive self-intersection.
No ALE or ALF K\"ahler 4-manifold can have a holomorphic, embedded $\mathbb{P}^1$ of positive self-intersection.

After this proposition and the two theorems, we lay out a few immediate consequence in the form of four corollaries.
A boundaryless K\"ahler manifold, compact or not, can have at most one embedded holomorphic $\mathbb{P}^1$ of positive self-intersection.
Any K\"ahler 4-manifold, compact or not, with an embedded holomorphic $\mathbb{P}^1$ of positive self-intersection has trivial fundamental group.
A classic result dating back to \cite{KR65} is that if every end of a K\"ahler manifold is ALE or ALF, then it has just one end.
We expand this classic result to show that such a manifold has only one end of ALE or ALF type, regardless of whatever other kinds of ends it might have, even if it has infinitely many ends.

The paper's final section has examples demonstrating the sharpness of some of these results.
This includes examples displaying several singular and non-singular Levi-flat submanifolds inside $\mathbb{P}^2$ and $\mathbb{P}^2\#\overline{\mathbb{P}}{}^2$---all of which, in the smooth case, intersect all curves of self-intersection $+1$.
Another example is an explicit construction of a complete, 2-ended scalar-flat K\"ahler metric, one of whose ends is ALE.
This is a K\"ahler metric conformally equivalent to the classic Taub-NUT metric.
Ordinarily in dimension 4 a K\"ahler metric cannot be conformal to another K\"ahler metric, but this is possible in the Taub-NUT case because of an unusually large number of metric-compatible complex structures available, including two that are not part of its hyperK\"ahler structure.
See Section \ref{SubSecTwoEnded}.

\subsection{Definitions and basic concepts}

The operator $\sqrt{-1}\partial\bar\partial:\bigwedge^0\rightarrow\bigwedge^2$ is real, and can be expressed $\sqrt{-1}\partial\bar\partial{}f=-\frac12d(J(df))$; throughout we use the convention $J\eta=\eta\circ{}J$ for 1-forms $\eta$ and we often use $dJdf$ for $d(J(df))$.
A function is called {\it pluriharmonic} when $\partial\bar\partial{}f=0$ or equivalently $dJdf=0$.

Let $V^3\subset{}M^4$ be a surface of class $C^2$, and let $\hat{n}$ be a unit 1-form defined in a neighborhood of $V^3$ that is perpendicular to $V^3$ in the sense that $\hat{n}(X)=0$ for every $X\in{}T_pV^3$.
Then the 2-form $L=-\frac12dJ\hat{n}\big|_{\hat{n}^\perp,J\hat{n}^\perp}\in\bigwedge^2V^3$ is called the Levi form of $V^{3}$.
Ostensibly $L$ depends a specific choice of $\hat{n}$ in a neighborhood of $L$, but it is easily shown to be independent of this choice, except that the sign of $L$ reverses if $\hat{n}$ is replaced by $-\hat{n}$.

The Levi form of a $C^2$ surface with respect to some normal is said to be non-negative, positive, or zero depending on whether $L(X,JX)\ge0$, $>0$, or $=0$, respectively, for all sections $X\in{}TV^3$ where also $JX\in{}TV^3$.
A surface $V^3$ with a choice of normal is defined to be \textit{pseudoconvex}, \textit{strictly pseudoconvex}, or \textit{Levi-flat} if $L\ge0$, $L>0$, or $L\equiv0$.
On a K\"ahler manifold, by dualizing we can evaluate $L$ on 1-forms rather than vectors.
Non-negativity means $L(J\eta,\eta)\ge0$ when $\eta,J\eta\perp\{\hat{n},J\hat{n}\}$---the reason for the switch in the place of $J$ from $L(X,JX)$ to $L(J\eta,\eta)$ is because of our convention for $J$: $J\eta$ means $+\eta\circ{}J$ instead of $-\eta\circ{}J$.
Given 1-forms $\eta,\gamma$ that are perpendicular to $\hat{n}$ and $J\hat{n}$ on $V^3$, $L(\eta,\gamma)$ is
\begin{eqnarray}
	L(\eta,\gamma)
	\;=\;-\frac12dJ\hat{n}\big(\eta,\gamma\big)
	\;=\;-\frac12*\big(\eta\wedge\gamma\wedge*dJ\hat{n}\big) \label{EqnNotationLWedge}
\end{eqnarray}
and we sometimes write $L(\eta\wedge\gamma)$ in place of $L(\eta,\gamma)$.

Given one or more boundary components, there are always bounded harmonic functions $f:M^4\rightarrow[0,1]$ that equal $0$ or $1$ on each boundary component.
Our main technical result is that when $f$ is such a function then
\begin{eqnarray}
	2\int_{M^4}|\partial\bar\partial{}f|^2
	\;=\;-\int_{\partial{}M^4}L\big(*(Jdf\wedge{}df)\big). \label{EqnIntByPartsIntro}
\end{eqnarray}
See (\ref{EqnCoerciveEst2}).
Clearly in the pseudoconvex case where $L\ge0$, we obtain some very coercive estimates---in particular harmonicity $\triangle{}f=0$ forces pluriharmonicity $\partial\bar\partial{}f=0$, and if $f$ is non-constant then $L\ge0$ forces $L\equiv0$.
Ordinarily harmonic functions are very easy to find and pluriharmonic functions are very hard to find.
The boundary condition $L\ge0$ makes the concepts much more closely tied together.
The following lemmas are almost immediate consequences of (\ref{EqnIntByPartsIntro}).
Unlike our theorems these two lemmas do not require $M^4$ to obey condition ($P$), which requires all boundary components to be compact.
Our lemmas only require some boundary components to be compact (but otherwise must obey ($P$)).
\begin{lemma} \label{LemmaHarmPluri}
	Assume $(M^4,g,J)$ satisfies ($P$) with the exception that its boundary components need not be compact, but that only a single boundary component $V^3\subseteq\partial{}M^4$ is compact.
	Let $f:M^4\rightarrow[0,1]$ be an harmonic function that equals $1$ on $V^3$ and $0$ on all other boundary components.
	
	Then $f$ is pluriharmonic.
	If $f$ is non-constant (eg. if $\partial{}M^4$ has more than one component, or $M^4$ is non-parabolic) then every component of $\partial{}M^4$ is Levi-flat.
\end{lemma}
\begin{lemma} \label{lemmaFundGrpsAndParaThm}
	Assume $(M^4,g,J)$ satisfies ($P$), and $V^3$ is a boundary component satisfying ($F$).
	
	Then $V^3$ is the {\it only} boundary component of $M^4$, $\iota_*(\pi_1(V^3))=\pi_1(M^4)$, and the Riemannian manifold $(M^4,g)$ is either parabolic or compact (in particular it has no ALE or ALF ends).
\end{lemma}
Lemmas \ref{LemmaHarmPluri} and \ref{lemmaFundGrpsAndParaThm} easily apply to the case of K\"ahler manifolds that have an ALE or ALF end.
Such manifold ends are always non-parabolic, because the Green's functions decay along such ends asymptotically like $dist^{-2}$ or $dist^{-1}$, respectively.
Such ends also meet the criteria from \cite{LT92} or \cite{HK01}, for example.
The paper \cite{Web19} contains results similar to these lemmas, arrived at differently.

\subsection{Results}

This paper focuses on how Levi-flat or pseudoconvex submanifolds within (or on the boundary of) a K\"ahler 4-manifold forces topological constraints on the ambient manifold.
This is analogous, in a certain sense, to the large body of literature on the subject of global curvature constraints forcing topological restrictions on manifolds.
The global condition in our setting is not global curvature controls, but rather the existence of a K\"ahler structure.
We use $\iota$ for the natural inclusion $\iota:V^3\rightarrow{}M^4$ of a submanifold into its ambient manifold.

\begin{theorem}[Levi-flat submanifolds] \label{ThmLeviFlats}
	Assume $(M^4,J,g)$ satisfies ($P$) and $V^3$ is a Levi-flat submanifold satisfying ($F$); in particular, $\pi_1(V^3)$ is finite.
	
	If $V^3$ is a boundary component then it is the only boundary component and $\iota_*\pi_1(V^3)=\pi_1(M^4)$.
	
	If $V^3$ is not a boundary component but an interior submanifold, then $M^4\setminus{}V^3$ has two components, $\partial{M}^4=\varnothing$, and $\iota_*\pi_1(V^3)=\pi_1(M^4)$.
	
	Finally, if $V_0^3$, $V_1^3$ are any two compact Levi-flat submanifolds that do not intersect $\partial{}M^4$ and both obey ($F$) except that only one but not both of $\pi_1(V_0^3)$, $\pi_1(V_1^3$) need be finite, then $V_0^3$ and $V_1^3$ intersect.
\end{theorem}

Our second result weakens the requirement that the submanifold be Levi-flat to the condition that it be pseudoconvex (with respect to either of its normals).
\begin{theorem}[Pseudoconvex submanifolds] \label{ThmPseudoconvexSubs}
	Assume $(M^4,J,g)$ satisfies ($P$), and $V^3$ is a pseudoconvex submanifold satisfying ($F$); in particular $\pi_1(V^3)$ is finite.
	
	If $V^3$ is a boundary component, then it is the only boundary component and $\iota_*\pi_1(V^3)=\pi_1(M^4)$.
	
	If $V^3$ is not a boundary component but an interior submanifold, then exactly one of the following holds:
	\begin{itemize}
		\item[{\it{i}})] $V^3$ is a separating surface; then letting $M_0^4\subset{}M^4\setminus{}V^3$ be the component with the pseudoconvex copy of $V^3$ on its boundary, we have $\partial{}M_0^4=V^3$ and $\iota_*\pi_1(V^3)=\pi_1(M^4_0)$.
		\item[{\it{ii}})] $V^3$ is a non-separating surface; then $\pi_1(M^4)=\iota_*\pi_1(V^3)\rtimes\mathbb{Z}$, $\partial{}M^4=\varnothing$, and the first betti number is unity: $b^1(M^4)=1$.
		In addition to $\partial{M}^4$ being boundaryless (and therefore geodesically complete), it is also non-compact.
	\end{itemize}
\end{theorem}
In particular, if a compact K\"ahler manifold has a pseudoconvex submanifold $V^3$ with finite fundamental group, $V^3$ is always separating.
Another immediate consequence is that any complete K\"ahler 4-manifold with one end that is ALE or ALF, it has finite fundamental group.

Theorem \ref{ThmLeviFlats} says, among other things, that non-intersecting Levi-flat submanifolds cannot exist in $M^4$ if either of their fundamental groups is finite.
If we weakening the assumption to pseudoconvexity then non-intersecting submanifolds certainly can exist, for instance concentric spheres in $\mathbb{C}^2$.
Nevertheless we can still draw some conclusions, which we summarize in the following proposition.
\begin{proposition}[Two pseudoconvex submanifolds] \label{PropEssSepPseudo}
	Let $(M^4,g,J)$ be a manifold (possibly non-compact) that satisfies ($P$).
	Assume $V_0^3$, $V_1^3$ are non-intersecting compact submanifolds that are pseudoconvex with respect to some orientation, and assume one of the fundamental groups, either $\pi_1(V^3_0)$ or $\pi_1(V^3_1)$, is finite.
	
	Then any component of $M^4\setminus\{V_0^3\cup{}V_1^3\}$ has, on its boundary, at most \emph{one} pseudoconvex copy of $V_0^3$ or $V_1^3$ but never both.
\end{proposition}
There is no apriori requirement that $V_0^3$ or $V_1^3$ be separating, no restrictions on the number of components $M^4\setminus\{V_0^3\cup{}V_1^3\}$ might have (it might have one, two, or three), and no restriction on the boundedness of the components.
It is also possible such a component may have pseudo-concave boundary components, which, by hypothesis ($P$), could only come from other copies of $V_0^3$ or $V_1^3$.

The proposition says nothing about, for instance, pairs of pseudoconvex submanifolds such as be concentric 3-spheres in $\mathbb{C}^2$ or perhaps small non-intersecting 3-spheres around distinct points in a K\"ahler manifold.
But as mild sounding as Proposition \ref{PropEssSepPseudo} may seem, it can have some strong consequences.
This theorem forbids any two ends of type ALE or ALF, regardless of what other kinds of ends the manifold might have---the reason is that ALE or ALF manifold ends have convex separating surfaces in the differential geometric sense (meaning they have positive second fundamantal forms) which are therefore also pseudoconvex.
(For a precise definition of ``ALE'' and ``ALF,'' see for example \cite{CK}.)
This theorem also gives extremely strong restrictions on K\"ahler 4-manifolds with an embedded, holomorphic $\mathbb{P}^1$ submanifold of \textit{positive} self-intersection, as any such $\mathbb{P}^1$ has a neighborhood which is pseudoconvex with respect to the \textit{inward} pointing normal.
Thus by Proposition \ref{PropEssSepPseudo} no two such $\mathbb{P}^1$ submanifolds can exist, unless they intersect.
If the K\"ahler manifold has an ALE or ALF end, no such $\mathbb{P}^1$ can exist at all.
This and other phenomena are explored in the Lemmas below.

\subsection{Consequences}

Our theorems produce certain topological restrictions on K\"ahler 4-manifolds from the existence of Levi-flat or pseudoconvex submanifolds.
Here we explore some immediate consequences of our theorems in four corollaries.
Much of the value in these corollaries lies in the fact that they apply in the compact as well as the non-compact case.
In the compact case the content of these lemmas is largely known.
For example, if a closed surface has a $\mathbb{P}^1$ of positive self-intersection it is already known to be a rational surface, and because these are classified, our Corollary \ref{CorCycsPseudo} for example was already known in the compact case.
\begin{corollary} \label{CorLeviFlatAndCycle}
	Assume $(M^4,J,g)$ is a complete K\"ahler manifold (not necessarily compact) that has both a compact embedded Levi-flat surface $V^3$ and an embedded holomorphic submanifold $N=\mathbb{P}^1$ of positive self-intersection.
	Then they intersect: $V^3\cap{}N\ne\varnothing$.
\end{corollary}
\begin{corollary} \label{CorCycsPseudo}
	Let $(M^4,J,g)$ be a pseudoconvex K\"ahler manifold (not necessarily compact) that satisfies ($P$).
	
	If $(M^4,J,g)$ has an embedded holomorphic $\mathbb{P}^1$ of positive self-intersection, then $\pi_1(M^4)=\{e\}$ and $\partial{}M^4=\varnothing$.
	
	If $(M^4,J,g)$ has an embedded holomorphic $\mathbb{P}^1$ of zero self-intersection, then every compact component of $\partial{}M^4$ has infinite $\pi_1(V^3)$.
\end{corollary}
\begin{corollary} \label{CorZeroAndPos}
	Assume $(M^4,J,g)$ is a complete K\"ahler manifold (not necessarily compact).
	If $N=\mathbb{P}^1$ and $N'=\mathbb{P}^1$ are embedded holomorphic submanifolds, one of which has non-negative self-intersection and one of which has positive self-intersection, then they intersect: $N\cap{}N'\ne\varnothing$.
\end{corollary}
\begin{corollary} \label{CorEnds}
	Let $(M^4,g,J)$ be a complete K\"ahler manifold that has $k$ many ALE ends and $l$ many ALF ends (apriori $k$ or $l$ might be infinite).
	Then $k+l$ is zero or one.
	If $k+l=1$ then $\pi_1(M^4)$ is finite.
	If $(M^4,g,J)$ has a holomorphic embedded $\mathbb{P}^1$ of self-intersection $\ge0$, then $k+l=0$.
\end{corollary}

An older theorem, originally a consequence of Kohn-Rossi's work \cite{KR65}, states that a K\"ahler manifold of finite topological type, all of whose ends are ALE or ALF, can have just one end.
More recently various other proofs of this theorem have emerged, using techniques very far from Kohn-Rossi's original methods.
In one work \cite{Web19} potential theory and the open mapping theorem was used.
In other works \cite{HLeB16}, \cite{LeB19} a compactification of ALE manifolds in the complex category was performed, which allowed the use of theorems designed for the compact case.
In still another \cite{LT92} potential theory was used, along with a curvature assumption used to control asymptotics of the Green's function.

\section{Interaction between harmonic functions and pseudoconvex boundary components} \label{SecHarmAndPseudos}

We shall frequently refer to harmonic functions created by the method of exhaustion; this is a common construction, and we run through it for the purposes of adaptation to our manifolds-with-boundary.
If $V_1^{3},\dots,V_K^{3}$ is a finite collection of compact boundary components we can create a bounded harmonic function $f$ with $f\equiv1$ on each $V_i^{3}$, $f\equiv0$ on all other components of $\partial{}M^4$.
To build this function, let $\{\Omega_i\}_{i=1}^\infty$ be pre-compact domains that each contain $\bigcup_{i=1}^KV_i^{3}$ and that exhaust $M^4$.
For each $i$ create harmonic functions $f_i:\Omega_i\rightarrow[0,1]$ with boundary conditions $f_i=1$ on $\bigcup_{i=1}^KV_i^{3}$ and $f_i=0$ on $\partial\Omega_i\setminus\bigcup_{i=1}^KV_i^{3}$.
Set $f_\infty=\lim_if_i$.
Then $f_\infty$ is certainly a bounded harmonic function $f_\infty:M^4\rightarrow[0,1]$ that is $1$ on $\bigcup_{i=1}^KV_i^{3}$ and zero on all other boundary components (if there are any others).
Such a function is unique, in the sense that if $f'_\infty$ is an harmonic function created in the same way except using a different exhaustion $\Omega_1'\subset\Omega_2'\subset\dots$, then $f'_\infty=f_\infty$.

Possibly an harmonic function created this way is constant.
This leads to the distinction between \textit{parabolic} and \textit{non-parabolic} manifold ends.
A non-compact Riemannian manifold with non-empty but compact boundary---a \textit{manifold end}---is called parabolic if the harmonic function created by the method of exhaustion is constant, and non-parabolic if it is non-constant.
This terminology, now standard, originated in \cite{LT92}.
If $M^n$ is a Riemannian manifold with non-empty but compact boundary and $f_\infty$ is the harmonic function created by this method of exhaustion, then the quantity $\int_{M^n}|df_\infty|^2$ is always finite, and is called the harmonic capacity of $M^n$.
This follows from a simple argument using the classical Hopf lemma: the pre-compact domains $\Omega_i$ all have finite capacity, and an application of the Hopf lemma shows that the value of $\int_{\Omega_i}|df_i|^2$ decreases as $i$ increases.
For more details, see for example \cite{LT92}.

\subsection{Integration by Parts}

We use the convention $\triangle{}f=+Tr\,\text{Hess}f$.
On a K\"ahler manifold with K\"ahler form $\omega$, the projection of $dJdf\in\bigwedge^2$ onto $\bigwedge^+$ is $-\frac12(\triangle{f})\omega$.
From this we obtain
\begin{eqnarray}
	*dJdf=-(\triangle{}f)\omega-dJdf. \label{EqnStarDJD}
\end{eqnarray}
Assume $M^4$ is a manifold with compact boundary, $f$ is a $C^2$ function, and $\varphi$ is a cutoff function---that is, a $C^\infty$ function with compact support on $\overline{M}{}^4$.
We do not require $\varphi$ to vanish on any particular boundary component.
Using (\ref{EqnStarDJD}), that $2\sqrt{-1}\partial\bar\partial{}f=-dJdf$, along with integration by parts we find
\begin{eqnarray}
	\begin{aligned}
		&4\int\varphi^2|\partial\bar\partial{}f|^2dVol
		\;=\;\int\varphi^2dJdf\wedge*dJdf \\
		&\quad
		\;=\;\frac12\int\varphi^2(\triangle{}f)^2\omega\wedge\omega
		-\int\varphi^2dJdf\wedge{}dJdf \\
		&\quad
		\;=\;\int\varphi^2(\triangle{}f)^2\,dVol
		+2\int\varphi\,d\varphi\wedge{}Jdf\wedge{}dJdf
		-\int_{\partial{}M^4}\varphi^2\,Jdf\wedge{}dJdf.
	\end{aligned} \label{EqnPrePreMainEst}
\end{eqnarray}
\begin{lemma}[Integration by Parts] \label{LemIBP}
	Let $(M^4,g)$ be a K\"ahler manifold with $C^2$ boundary.
	Let $V_1^3,\dots,V_K^3\subset\partial{}M^4$ be finitely many components of $\partial{}M^4$, each of which is compact (although $\partial{}M^4$ might have other components, possibly infinitely many, that might not be compact).
	Let $f:M^4\rightarrow[0,1]$ be the harmonic function obtained from the exhaustion method, where $f\equiv1$ on each $V_1^3,\dots,V_K^3$ and $f\equiv0$ on every other boundary component.	
	Let $\varphi$ be any $C_c^{0,1}$ function with compact support on $\overline{M}{}^4$.
	Then
	\begin{equation} \label{EqnCoerciveEst1}
		2\int\varphi^2|\partial\bar\partial{}f|^2
		\;=\;\int\varphi\,d\varphi\wedge{}Jdf\wedge{}dJdf
		-\int_{\partial{}M^4}\varphi\,L\big(*(Jdf\wedge{}df)\big)\,dA
	\end{equation}
	where $L$ is the Levi form of $\partial{}M^4$.
\end{lemma}
{\it Proof.}
When $f:M^4\rightarrow[0,1]$ is harmonic and $\varphi$ has compact support, (\ref{EqnPrePreMainEst}) is
\begin{eqnarray}
	\begin{aligned}
		&4\int\varphi^2|\partial\bar\partial{}f|^2
		\;=\;2\int\varphi\,d\varphi\wedge{}Jdf\wedge{}dJdf
		-\int_{\partial{}M^4}\varphi\,Jdf\wedge{}dJdf.
	\end{aligned} \label{EqnPreMainEst}
\end{eqnarray}
We relate the boundary term of (\ref{EqnPreMainEst}) to the Levi form.
First, we rewrite
\begin{eqnarray}
	-\int_{\partial{}M}\varphi\,Jdf\wedge{}dJdf
	\;=\;
	-\int_{\partial{}M}\varphi
	*\big(
	\hat{n}\wedge{}Jdf\wedge{}dJdf
	\big)\,dA \label{EqnBoundTermFirst}
\end{eqnarray}
where $dA$ is the area form on $\partial{}M^4$, $\hat{n}$ is the outward unit normal, and $*$ is the Hodge-star on $M^4$.
Because $f:M^4\rightarrow[0,1]$ is harmonic and either zero or one on each boundary component, certainly $|df|>0$ and $f$ is a defining function for each component.
Therefore the Levi form is
\begin{eqnarray}
	L
	\;=\;\mp\frac12dJ\left(|df|^{-1}df\right)\Big|_{df^\perp,Jdf^\perp}
	\;=\;\mp\frac12|df|^{-1}\,dJdf\Big|_{df^\perp,Jdf^\perp}, \label{EqnLeviBoundary}
\end{eqnarray}
where $\mp$ is $-$ when $df$ is outward pointing (this is the case when $f=1$ on a boundary component, which is the case on all the $V_i^3$ components), and is $+$ when $df$ is inward pointing (this is the case when $f=0$ on a boundary component, which is the case on all components of $\partial{}M^4$ except the $V_i^3$).

The outward pointing 1-form is $\hat{n}=\pm|df|^{-1}df$.
Using the fact that $*\omega_1\wedge*\omega_2=\omega_1\wedge\omega_2$ whenever $\omega_1,\omega_2\in\bigwedge{}^2$, the integrand of (\ref{EqnBoundTermFirst}) is
\begin{eqnarray}
	\begin{aligned}
		\hat{n}\wedge{}Jdf\wedge{}dJdf
		&\;=\;\pm\,df\wedge{}Jdf\wedge\left(|df|^{-1}dJdf\right) \\
		&\;=\;\pm*(df\wedge{}Jdf)\wedge*\left(|df|^{-1}dJdf\right) \\
		&\;=\;-2*\,L\big(*(df\wedge{}Jdf)\big) \\
		&\;=\;+2*\,L\big(*(Jdf\wedge{}df)\big).
	\end{aligned} \label{EqnEvaluatingLevi}
\end{eqnarray}
(The ``$\pm$'' from $\hat{n}=\pm|df|^{-1}df$ and the ``$\mp$'' of (\ref{EqnLeviBoundary}) combine to make ``$-$''.)
Therefore (\ref{EqnPreMainEst}) is precisely
\begin{eqnarray}
	2\int\varphi^2|\partial\bar\partial{}f|^2
	\;=\;\int\varphi\,d\varphi\wedge{}Jdf\wedge{}dJdf
	-\int_{\partial{}M^4}\varphi\,L\big(*(Jdf\wedge{}df)\big).
\end{eqnarray}
\qed

\begin{lemma}[Improved integration by parts] \label{LemmaImprIBP}
	Let $(M^4,g)$ be a K\"ahler manifold along with an harmonic function $f:M^4\rightarrow[0,1]$ that satisfy the hypotheses of Lemma \ref{LemIBP}.
	Then if $\int_{\partial{}M^4}L\big(*(Jdf\wedge{}df)\big)\,dA$ is integrable, 
	\begin{equation} \label{EqnCoerciveEst2}
		2\int_{M^4}|\partial\bar\partial{}f|^2
		\;=\;-\int_{\partial{}M^4}L\big(*(Jdf\wedge{}df)\big)\,dA.
	\end{equation}
	Finally, if $\partial{}M^4$ is pseudoconvex then $\int_{\partial{}M^4}L\big(*(Jdf\wedge{}df)\big)\,dA$ is integrable.
\end{lemma}
{\it Proof.}
We prove this using (\ref{EqnCoerciveEst1}) and the fact that the capacity $\int|df|^2$ is finite.
The method is by choosing good cutoff functions.
We are brief because this kind of argument is very standard.
For each $i$ let $\varphi_i$ be a cutoff function with $\varphi_i=1$ on a very large compact set $\Omega_i$, $\varphi=0$ on an even larger compact set $\Omega_{i+1}$, and the exhaustion $\Omega_1\subset\Omega_2\subset\dots$ is chosen so also $|d\varphi_i|<i^{-1}$.
It might be objected that our manifolds have boundary, so the exponential map cannot provide the distance functions necessary to build such cutoff functions $\varphi_i$.
However a manifold with $C^2$ boundary still has a distance function to any point, given by approximation by minimizing paths.
At interior points of the manifold such a distance function retains the usual properties, in particular being Lipschitz and having unit norm almost everywhere.
This is enough to construct cutoff functions in the usual way.

We estimate the first integral on the right of (\ref{EqnCoerciveEst1}), using H\"older's inequality, by
\begin{eqnarray}
	\begin{aligned}
		&\left|\int_{M^4}\varphi_i{}d\varphi_i\wedge{}Jdf\wedge{}dJdf\right|
		\;\le\;
		\max|d\varphi_i|\int_{M^4}\varphi_i|df||dJdf| \\
		&\hspace{1in}
		\;\le\;\frac12\max|d\varphi_i|\left(\int_{M^4}|df|^2
		\,+\,\int_{M^4}\varphi_i^2|dJdf|^2\right).
	\end{aligned}
\end{eqnarray}
Because $|d\varphi_i|<i^{-1}$, we have upper and lower estimates
\begin{eqnarray}
	\begin{aligned}
		&2\left(1-i^{-1}\right)\int\varphi^2_i|\partial\bar\partial{}f|^2
		\;\le\;\frac12i^{-1}\int|df|^2
		\,-\,\int_{\partial{}M^4}\varphi_i{}L\Big(*(Jdf\wedge{}df)\Big)\,dA \\
		&2\left(1+i^{-1}\right)\int\varphi^2_i|\partial\bar\partial{}f|^2
		\;\ge\;-\frac12i^{-1}\int|df|^2
		\,-\,\int_{\partial{}M^4}\varphi_i{}L\Big(*(Jdf\wedge{}df)\Big)\,dA.
	\end{aligned}
\end{eqnarray}
But the capacity term $\int|df|^2$ is finite.
Therefore, as long as $\int_{\partial{}M^4}L\big(*(Jdf\wedge{}df)\big)\,dA$ is integrable, we take $i\rightarrow\infty$ and obtain (\ref{EqnCoerciveEst2}).

To see this is integrable when $\partial{}M^4$ is pseudoconvex, recall from the introduction that this means $L(J\eta,\eta)\ge0$ when $\eta,J\eta$ are perpendicular to $\hat{n}$ and $J\hat{n}$.
Thus, because $*(Jdf\wedge{}df)=J\eta\wedge\eta$ for some 1-form $\eta$, then $L\big(*(Jdf\wedge{}df)\big)\ge0$.
Because the integrand is pointwise non-negative, $\int_{\partial{}M^4}\varphi_i{}L\big(*(Jdf\wedge{}df)\big)\,dA$ is integrable (although it might equal $+\infty$).
This concludes the proof.
\qed

Compare the argument of Lemma \ref{LemmaImprIBP} to the argument of Lemma 3.1 of \cite{Li90}.

\textit{Proof of Lemma \ref{LemmaHarmPluri}.}
Our assumption is that $\partial{}M^4$ is entirely pseudoconvex, and $V^3\subseteq\partial{}M^4$ is a compact boundary component.
Create an harmonic function $f$ by the method of exhaustion so $f$ is unity on $V^3$ and zero on any other components of $\partial{}M^4$.
Possibly $f$ is constant, in which case $f\equiv1$ and $df\equiv0$.
Certainly in this case $\partial{M}^4$ has a single boundary component, and is either compact or parabolic.

Possibly $f$ is non-constant.
Of so, then by pseudoconvexity, $L(*(Jdf\wedge{}df))\ge0$ in the pointwise sense on $\partial{}M^4$.
Therefore (\ref{EqnCoerciveEst2}) forces both $|\partial\bar\partial{}f|\equiv0$ (so $f$ is pluriharmonic), and $L(*(df\wedge{}Jdf))\equiv0$.
By the Hopf lemma $df$ is never zero on $\partial{}M^4$ because every point on $\partial{}M^4$ is either a global maximum or a global minimum of $f$.
Thus $L(*(df\wedge{}Jdf))\equiv0$ forces $L\equiv0$ on $\partial{}M{}^4$, which is the same as $\partial{}M{}^4$ being Levi-flat.
This concludes the proof of Lemma \ref{LemmaHarmPluri}.

\begin{lemma}[Analytic maps into the strip] \label{LemmaAnaStrip}
	Let $(M^4,J)$ be a compact complex manifold with precisely two boundary components $V_0^3$ and $V_1^3$.
	If $C_0\ne{}C_1$ there is no analytic function $z:M^4\rightarrow\mathbb{C}$ with $Re(z)\big|_{V_0^3}\equiv{}C_0$ and $Re(z)\big|_{V_1^3}\equiv{}C_1$.
\end{lemma}
\textcolor{white}{ }\indent \textbf{Remark}.
By ``strip'' we mean a locus $\{z\in\mathbb{C}\,|\,C_0\le{}Re(z)\le{}C_1\}$.

{\it Proof.}
This simple lemma follows from the open mapping theorem.
With $z:M^4\rightarrow\mathbb{C}$ being such an analytic function, by continuity the image $z(M^4)$ is connected and compact.
By the open mapping theorem, $\partial(z(M^4))\subseteq{}z(\partial{}M^4)$.

But by hypothesis, the image of the boundary $z(\partial{}M^4)$ is contained within the the lines $\{Re(z)=C_0\}$ and $\{Re(z)=C_1\}$, and intersects \textit{both} of these lines.
Because $z(M^4)$ is connected and intersects both lines, it intersects interior points of the strip.
Because $z(M^4)$ is closed it is either the entire strip or else it has points of closure in the interior of the strip.
Because $z(M^4)$ is compact it cannot be the entire strip, so $z(M^4)$ must have points of closure within the strip.
This contradicts the fact that the image's boundary $\partial(z(M^4))\subseteq{}z(\partial{}M^4)$ lies completely within the strip's boundary lines, and concludes the proof.
\qed

{\it Proof of Lemma \ref{lemmaFundGrpsAndParaThm}.}
The hypotheses are that $(M^4,J,g)$ is K\"ahler with pseudoconvex boundary, and has a boundary component $V^3\subseteq\partial{}M^4$ that is compact and has finite $\pi_1(V^3)$.

We first show that $\iota_*(\pi_1(V^3))=\pi_1(M^4)$.
Let $\widetilde{M}{}^4$ be the universal cover of $M^4$ with covering map $p:\widetilde{M}{}^4\rightarrow{}M^4$.
Let $\widetilde{V}^3=p^{-1}(V^3)$ be the pre-image of $V^3$.
Each component of $\widetilde{V}^3$ covers $V^3$ with deck group $\iota_*(\pi_1(V^3_i))$, and because $\iota_*(\pi_1(V^3))$ is finite each component is itself compact.
The number of components is the cardinality of the coset space $\pi_1(M^4)/\iota_*(\pi_1(V^3))$.

For a proof by contradiction assume the inclusion $\iota_*(\pi_1(V^3))\subseteq\pi_1(M^4)$ is strict, so the pre-image $\widetilde{V}{}^3=p^{-1}(V^3)$ has at least two components which we label $\widetilde{V}{}_0^3,\widetilde{V}{}_1^3\subset{}\widetilde{V}{}^3$.
Let $f:\widetilde{M}{}^4\rightarrow[0,1]$ be the harmonic function created by the method of exhaustion, with $f\equiv1$ on $\widetilde{V}{}_1^3$ and $f\equiv0$ on all other components of $\partial\widetilde{M}{}^4$ including $\widetilde{V}{}_0^3$.
By hypothesis the boundary integral on the right side of (\ref{EqnCoerciveEst2}) is non-negative (indeed, even the integrand is non-negative in the pointwise sense), so the equality
\begin{equation}
	2\int_{\widetilde{M}{}^4}|\partial\bar\partial{}f|^2
	\;=\;-\int_{\partial\widetilde{M}{}^4}L\big(*(Jdf\wedge{}df)\big)\,dA
\end{equation}
now forces both sides to equal zero.
Because $f$ is not constant (it equals $1$ on $\widetilde{V}_1^3$ and $0$ on $\widetilde{V}_0^3$), we have $Jdf\wedge{}df$ nowhere zero on the boundary.
This forces $L\equiv0$ on $\partial{}M^4$, and also forces $\partial\bar\partial{}f=0$ on $M^4$ which means $f$ is pluriharmonic.
Because $dJdf=0$ (which is the same as $\partial\bar\partial{}f=0$) and because $\widetilde{M}{}^4$ is simply connected, the equation $dJdf=0$ means $Jdf=-dg$ for some real-valued function $g$.
The function $z=f+\sqrt{-1}g$ is therefore an analytic function.

Consider the submanifold
\begin{eqnarray}
	M_\epsilon^4
	\;=\;\left\{p\in\widetilde{M}{}^4\;\big|\;
	f(p)\in[1-\epsilon,1] \right\}.
\end{eqnarray}
Because $f=1$ and $df\ne0$ on $\widetilde{V}_1^3$ and $\widetilde{V}_1^3$ this is a global maximum for $f$, we can choose $\epsilon>0$ so small that $M_\epsilon^4$ is a collar neighborhood: $M_\epsilon^4\approx\widetilde{V}_1^3\times[1-\epsilon,1]$.
But then $M_\epsilon^4$ is a manifold with two compact boundary components, and the analytic function $z:M_\epsilon^4\rightarrow\mathbb{C}$ has $Re(z)=f=1-\epsilon$ on one boundary and $Re(z)=f=1$ on the other boundary.
By Lemma \ref{LemmaAnaStrip} applied to $\widetilde{M}{}^4_\epsilon$ along with the analytic function $z$, this is impossible.
We conclude that $f:\widetilde{M}{}^4\rightarrow[0,1]$ must be a constant function.
Thus two components $\widetilde{V}_0^3$ and $\widetilde{V}_1^3$ of $p^{-1}(V^3)$ cannot exist.
Therefore $\iota_*(\pi_1(V^4))=\pi_1(M^4)$.

Next we show that actually $\partial{}M^4=V^3$.
Let $p:\widetilde{M}{}^4\rightarrow{}M^4$ again be the universal cover.
We have shown $\iota_*(\pi_1(V^3))=\pi_1(M^4)$, so $p$ is a finite-sheeted cover.
By assumption the components of $\partial{}M^4$ are compact, so all components of $\partial\widetilde{M}{}^4$ remain compact.
Repeating the argument above, we conclude that $\partial\widetilde{M}{}^4$ has just one component, $p^{-1}(V^3)$.
Therefore $V^3$ is the only component of $\partial{}M^4$.

Lastly we must show that $(M^4,J,g)$ is parabolic or compact.
The boundary $\partial{}M^4$ has just one component, which is compact and pseudoconvex.
If $M^4$ is non-parabolic, then the harmonic function $f:M^4\rightarrow[0,1]$ created by the method of exhaustion is non-constant and a defining function on $\partial{}M^4$.
Since the boundary $\partial{}M^4$ is compact, we can repeat the previous argument---since $f$ is non-constant and reaches a maximum at $\partial{}M^4$, we can still use a collar neighborhood of $\partial{}M^4$ and an analytic function $z$ with $Re(z)=f$, and draw a contradiction using Lemma \ref{LemmaAnaStrip}.
Therefore $M^4$ cannot be non-parabolic, so it is parabolic or compact.
This concludes the proof of Lemma \ref{lemmaFundGrpsAndParaThm}.

\section{Levi-flat and Pseudoconvex Submanifolds}

Up to this point, our results have dealt with boundary pseudoconvexity (see \cite{Web19} for similar results using different methods).
We now move from boundary pseudoconvexity to pseudoconvex submanifolds.
To study submanifolds, we cut $M^4$ along these submanifolds, which produces additional boundary components, and using these new boundary components it is often possible to apply our lemmas.

{\it Proof of Theorem \ref{ThmLeviFlats}.}
The hypotheses are that $(M^4,J,g)$ is a K\"ahler manifold with pseudoconvex (or empty) boundary, has a compact embedded Levi-flat submanifold $V^3$ with finite $\pi_1(V^3)$, and either $V^3$ is a boundary component or else $V^3\cap\partial{}M^4=\varnothing$.

If $V^3$ is a boundary component, then the fact that $\iota_*(\pi_1(V^3))=\pi_1(M^4)$ follows from Lemma \ref{lemmaFundGrpsAndParaThm}.
The rest of the proof deals with the case $V^3\cap\partial{}M^4=\varnothing$.

\underline{\it Proof that $V^3$ is separating}.
For a contradiction, assume $V^3$ is non-separating.
This means $M^4{}'=M^4\setminus{}V^3$ consists of a just one component.
Because $V^3$ and $M^4$ are orientable (so in particular the submanifold $V^3$ has two distinct sides), the connected manifold $M^4{}'$ has two additional boundary components, each a copy of $V^3$, which we call $V^3_0$ and $V^3_1$.
Because the submanifold $V^3$ was Levi-flat, the new boundary components $V_0^3$ and $V_1^3$ are also Levi-flat.
Therefore the boundary of $M^4{}'$ continues to be pseudoconvex, so $M^4{}'$ continues to obey condition ($P$).
Both $V^3_0$ and $V^3_1$ have finite fundamental groups (as both are diffeomorphic to $V^3$).
This contradicts Lemma \ref{lemmaFundGrpsAndParaThm}, which says $\partial{}M^4{}'$ can have at most one such component.
This contradiction forces the submanifold $V^3$ to be separating.

\underline{\it Proof that $\partial{}M^4=\varnothing$.}
Because $V^3$ is separating, $M^4\setminus{}V^3$ consists of two components which we call $M_0^4{}'$ and $M_1^4{}'$.
In addition to inheriting any boundary components of the original manifold $M^4$, the new manifolds-with-boundary $M_0^4{}'$, $M_1^4{}'$ each has as an additional boundary component which is a copy of $V^3$, which we call $V_0^3\subseteq\partial{}M^4_0{}'$ and $V_1^3\subseteq\partial{}M^4_1{}'$.
The components $M_i^4{}'$ continue to satisfy ($P$), because the new boundary components are Levi-flat.
But by Lemma \ref{lemmaFundGrpsAndParaThm}, $V_i^3$ is the only boundary component of $M_i^4{}'$, $i=1,2$.
Because each $M_i^4{}'$ has no boundary components besides $V_i^3$, the original manifold $M^4$ can have no boundary components at all.

\underline{\it Case that $V^3_0$, $V^3_1$ are non-intersecting.}
Specifically, the assumption is that $V^3_0$, $V^3_1$ are non-intersecting and Levi-flat, neither intersects $\partial{}M^4$, and at least one of them has finite fundamental group.

Without loss of generality we may assume $V^3_0$ has finite fundamental group.
From above, we know $\partial{}M^4=\varnothing$ and we know $V^3_0$ is separating, although we do not know if $V^3_1$ is separating or not.
Then $M^4\setminus(V_0^3\cup{}V_1^3)$ has either has three components (if $V_1^3$ is separating) or two components (if $V_1^3$ is not separating).
In either case, there is one component $M^4{}'$ of $M^4\setminus(V_0^3\cup{}V_1^3)$ that has on its boundary both a copy of $V_0^3$ and also a copy of $V_1^3$ (if $V_1^3$ is not separating, then $M^4{}'$ actually has two copies of $V_1^3$ on its boundary).

Because both $V_0^3$ and $V_1^3$ are Levi-flat, the component $M^4{}'$ continues to have pseudoconvex boundary.
But because the copy of $V_0^3$ on its boundary has finite fundamental group, Lemma \ref{lemmaFundGrpsAndParaThm} states this is the only boundary component, contradicting the existence of $V_1^3$.

\underline{\it Proof that $\iota_*(\pi_1(V^3))=\pi_1(M^4)$.}
The assumption is that $V^3$ is a Levi-flat submanifold that does not intersect $\partial{}M^4$, and that $\pi_1(V^3)$ is finite.

For a proof by contradiction suppose the inclusion $\iota_*(\pi_1(V^3))\subseteq\pi_1(M^4)$ is strict.
Because this inclusion is strict, after passing to the universal cover $p:\widetilde{M}{}^4\rightarrow{}M^4$ the preimage $\widetilde{V}{}^3=p^{-1}(V^3)$ has more than one component.
Let $\widetilde{V}{}^3_0,\widetilde{V}{}^3_1\subset{}p^{-1}(V^3)$ be two distinct components.
Each component of $p^{-1}(V^3)$ is itself a cover of $V^3$, and so the $V_i^3$ both have finite fundamental group.
But the result above states that no two Levi-flat submanifolds can exist in $M^4$ unless they intersect; this contradicts the fact that $\widetilde{V}{}^3_0$ and $\widetilde{V}{}^3_1$ are distinct components of $p^{-1}(V^3)$.
We conclude that $\iota_*\pi_1(V^3)=\pi_1(M^4)$, and conclude the proof of Theorem \ref{ThmLeviFlats}.

\vspace{0.1in}

\noindent{\it Proof of Proposition \ref{PropEssSepPseudo}.}
By hypothesis we have two non-intersecting pseudoconvex submanifolds $V_0^3$ and $V_1^3$.
We must show that if $M^4{}'$ is any component of $M^4\setminus(V_0^3\cup{}V_1^3)$, then $\partial{}M^4{}'$ contains has at most one pseudoconvex copy of $V_0^3$ or $V_1^3$.

\begin{figure}[h]
	\begin{subfigure}[b]{0.29\textwidth}
		\centering
		\includegraphics[width=\textwidth]{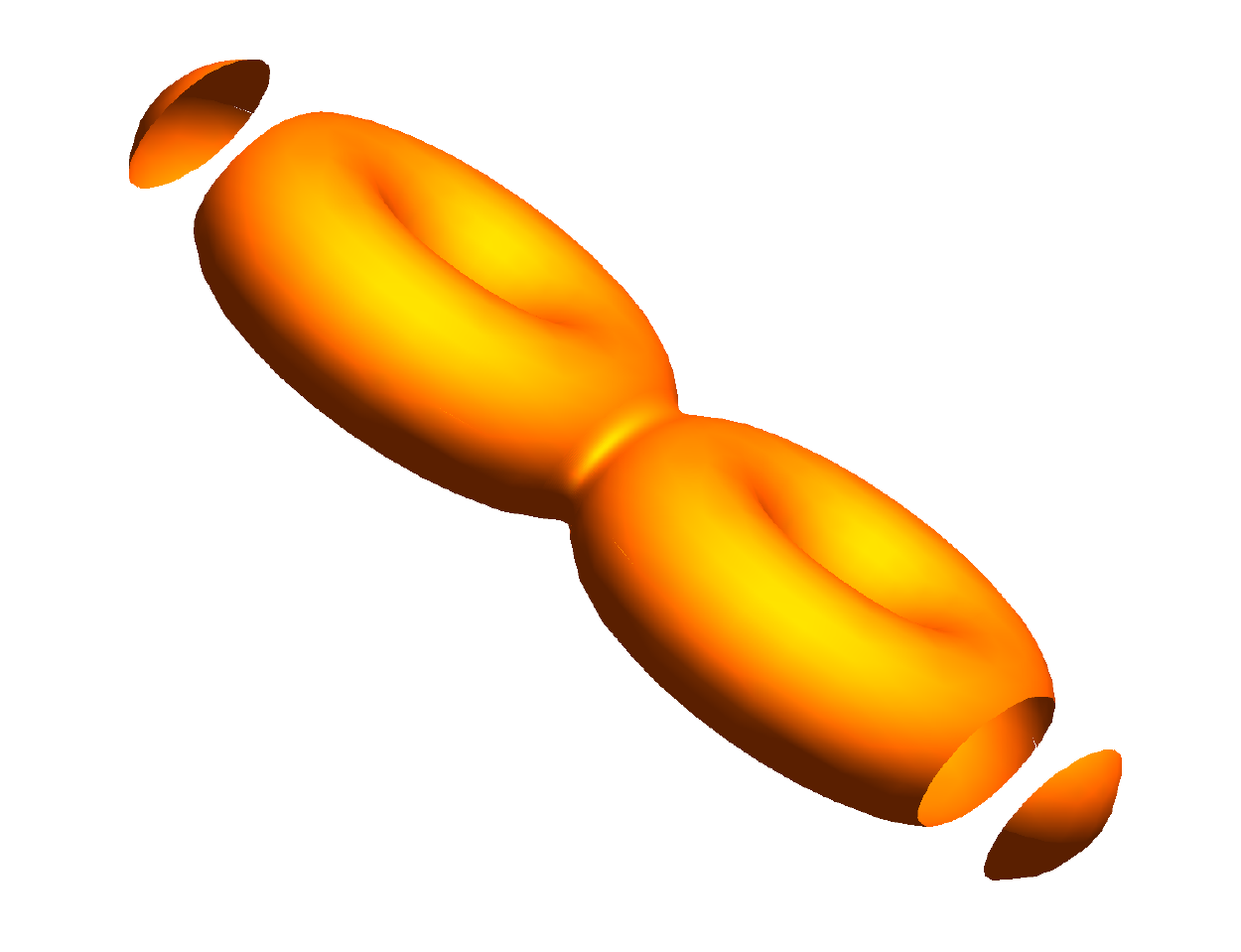}
		\caption{Case ({\it{i}}): $\partial{}M^4{}'$ has one pseudoconvex copy each of $V_0^3$ and $V_1^3$, and no pseudo-concave copies.
		}
	\end{subfigure}
	\hfill
	\begin{subfigure}[b]{0.29\textwidth}
		\centering
		\includegraphics[width=\textwidth]{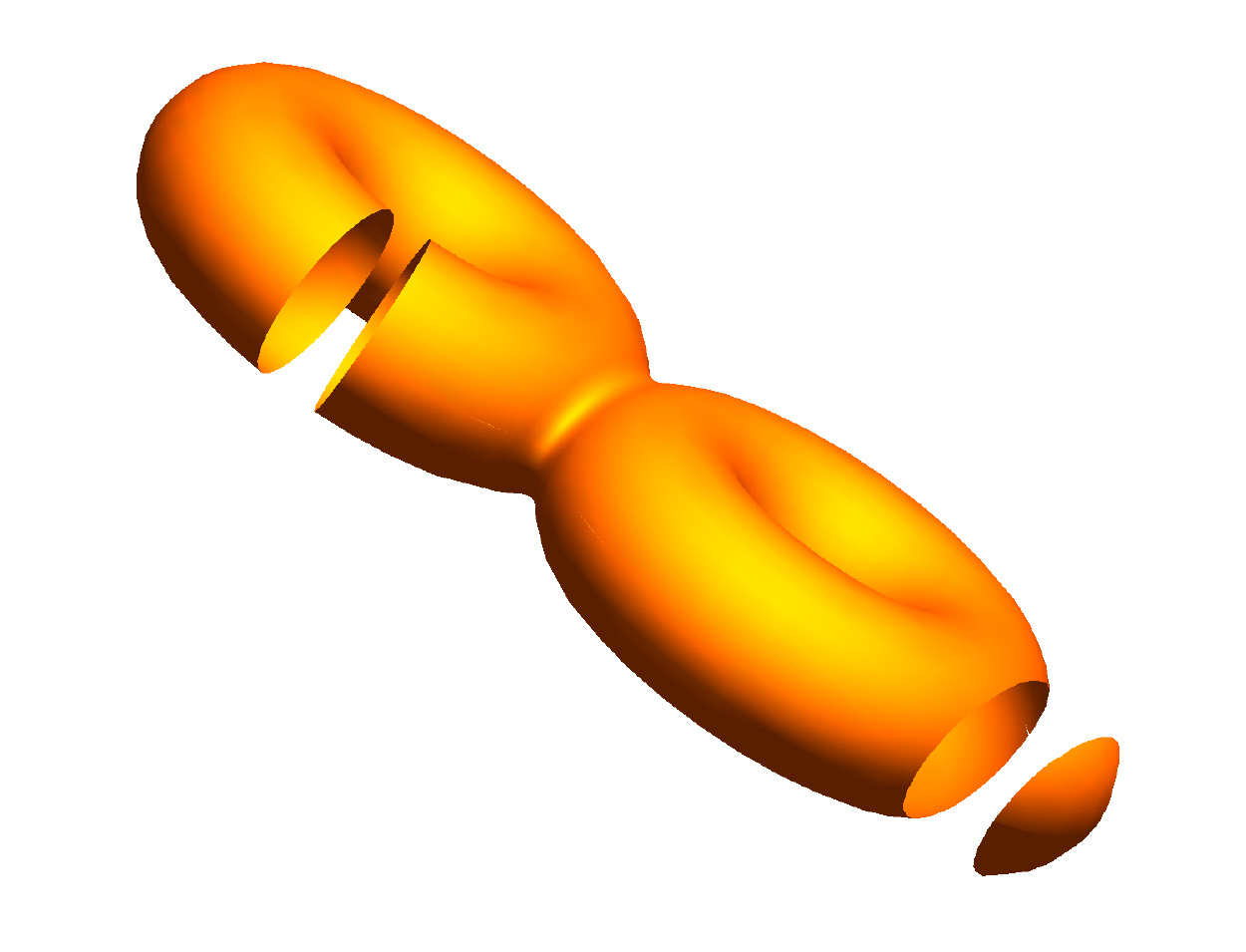}
		\caption{Case ({\it{ii}}): $\partial{}M^4{}'$ has one pseudoconvex copy of both $V_0^3$ and $V_1^3$, and one pseudo-concave copy of $V_0^3$.}
	\end{subfigure}
	\hfill
	\begin{subfigure}[b]{0.29\textwidth}
		\centering
		\includegraphics[width=\textwidth]{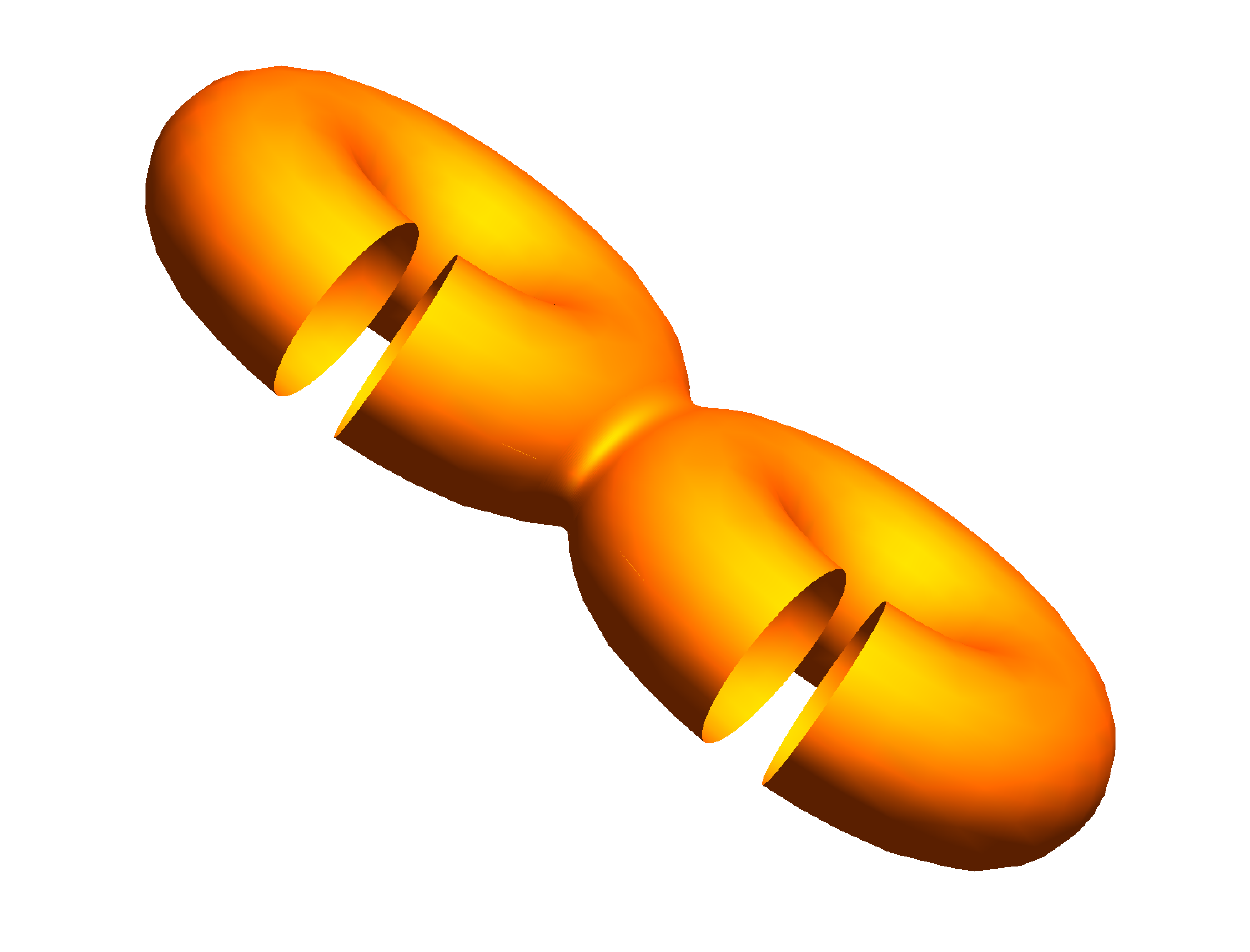}
		\caption{Case ({\it{iii}}): $\partial{}M^4{}'$ has one pseudoconvex and one pseudo-concave copy each of $V_0^3$ and $V_1^3$.}
	\end{subfigure}
	\caption{Schematic depiction of possibilities ({\it{i}}), ({\it{ii}}), ({\it{iii}}) for $M^4\setminus(V_0^3\cup{}V_1^3)$.}
	\label{FigThreeCases}
\end{figure}

Assuming not.
Then a component $M^4{}'$ of $M^4\setminus(V_0^3\cup{}V_1^3)$ has copies of both $V_0^3$ and $V_1^2$ on its boundary so that one copy of $V_0^3$ and one copy of $V_1^3$ is pseudoconvex with respect to the outward normal.
There are exactly three possibilities:
\begin{itemize}
	\item[{\it{i}})] $\partial{}M^4{}'$ has one pseudoconvex copy each of $V_0^3$ and $V_1^3$, and no pseudo-concave components.
	\item[{\it{ii}})] $\partial{}M^4{}'$ has one pseudoconvex copy each of $V_0^3$ and $V_1^3$ and one pseudo-concave copy of either $V_0^3$ or $V_1^3$ but not both.
	\item[{\it{iii}})] $\partial{}M^4{}'$ has one pseudoconvex and one pseudo-concave copy each of $V_0^3$, $V_1^3$.
\end{itemize}
See Figure \ref{FigThreeCases} for a depiction.
Of course, $\partial{}M^4{}'$ may contain additional pseudoconvex components, inherited from the original manifold $M^4$.

\underline{\it Proof in case ({\it{i}}).}
In this case since $M^4{}'$ has an entirely pseudoconvex boundary, with compact boundary components.
The fact that either $\pi_1(V_0^3)$ or $\pi_1(V_1^3)$ is finite contradicts Lemma \ref{lemmaFundGrpsAndParaThm}.

\underline{\it Construction of the pseudoconvex manifold $M_\infty^4$}.
In cases ({\it{ii}}) and ({\it{iii}}) the manifold $M^4{}'$ does not have pseudoconvex boundary, because $\partial{}M^4$ contains pseudo-concave copies of $V_0^3$ and/or $V_1^3$.

\begin{figure}[h]
	\vspace{-0.3in}
	\begin{tabular}{ll}
		\includegraphics[width=0.45\textwidth]{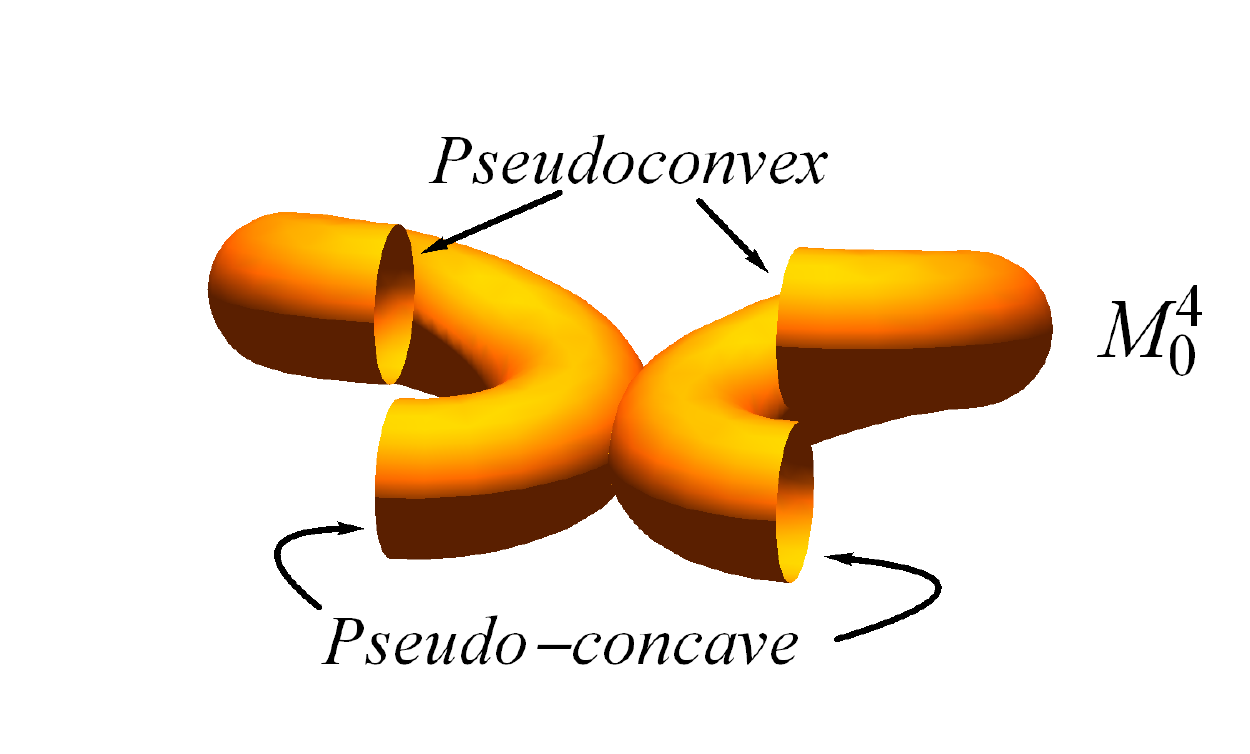} &
		\raisebox{2.4\height}{\parbox{0.5\textwidth}{
				(a) In cases ({\it{ii}}) and ({\it{iii}}) the boundary $\partial{}M_0^4$ has both pseudo-concave and pseudoconvex copies of $V_0^3$ and/or $V_1^3$.}} \\
		\includegraphics[width=0.45\textwidth]{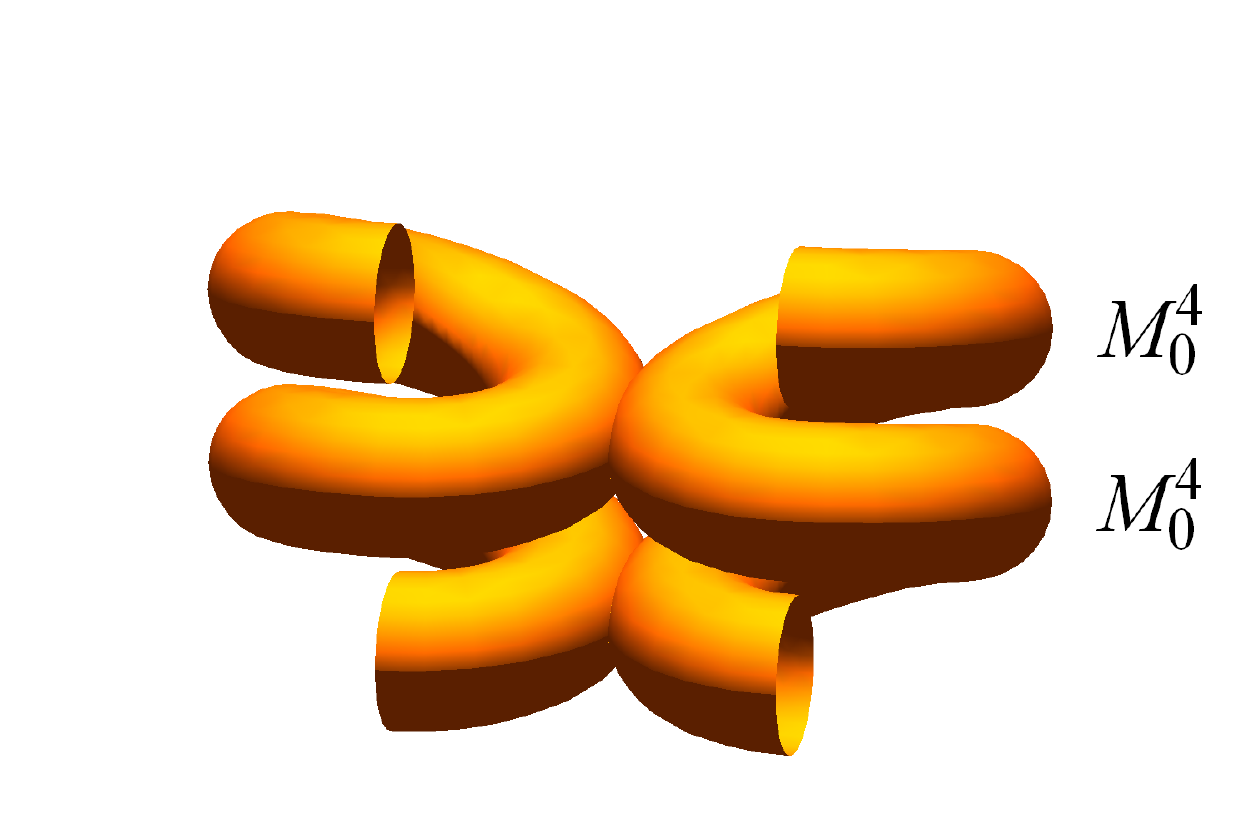} &
		\raisebox{1.0\height}{\parbox{0.5\textwidth}{
				(b) The manifold-with-boundary $M_1^4$ is created by gluing the pseudo-concave boundary components of $M_0^4$ to the matching pseudoconvex boundary components on another copy of $M_0^4$.
				Certainly $M_1^4=M_0^4\cup{}M_0^4$ continues to have both pseudo-concave and pseudo-convex boundary components.
		}} \\
		\includegraphics[width=0.45\textwidth]{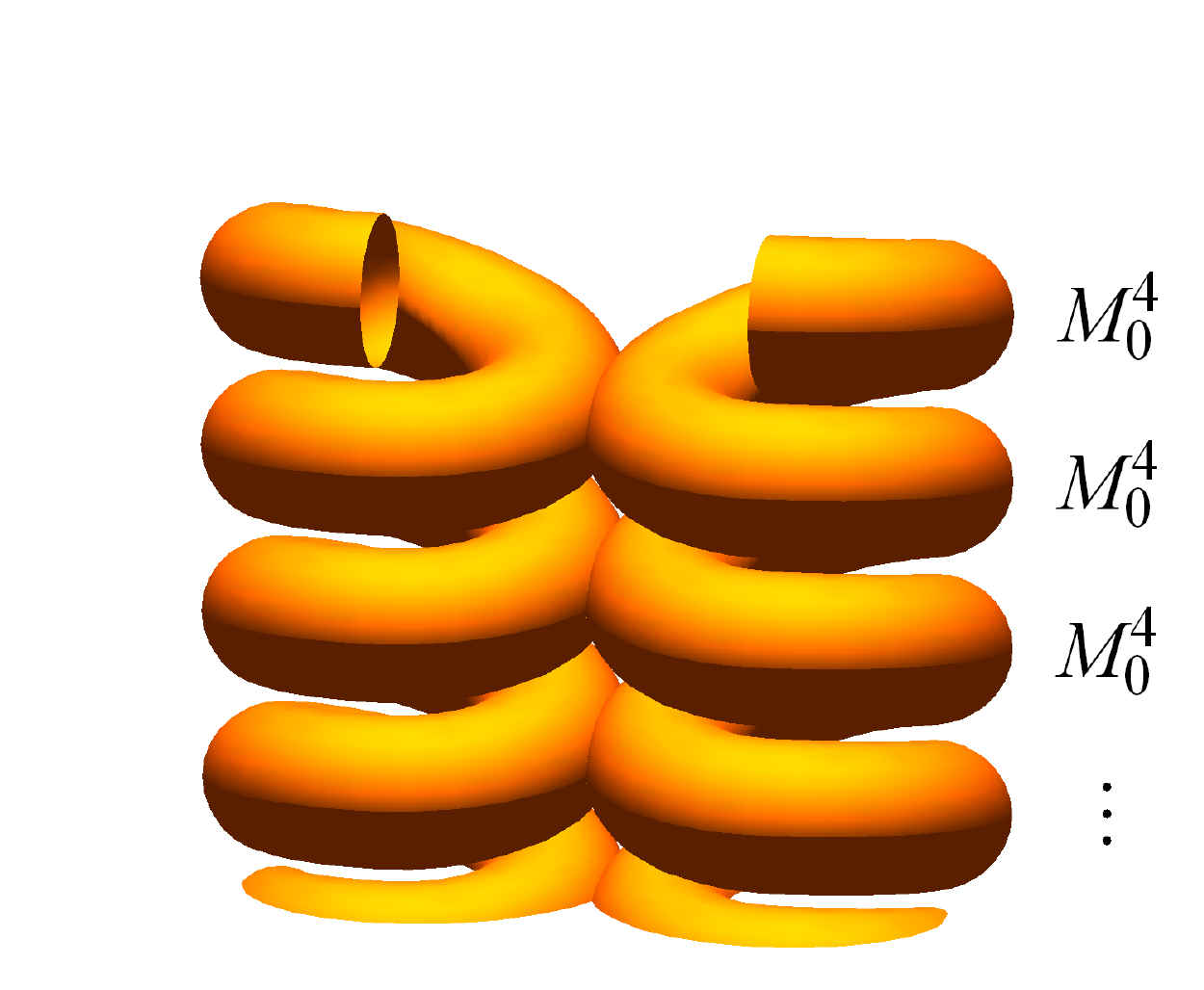} &
		\raisebox{1.4\height}{\parbox{0.5\textwidth}{
				(c) Having formed $M^4_{i-1}$, we glue the pseudo-concave components of $\partial{}M_{i-1}^4$ to the matching pseudoconvex components of yet another $M_0^4$ to create $M_i^4=M_0^4\cup\dots\cup{}M_0^4$ (union of $i$ many copies of $M_0^4$).
				Continuing indefinitely, we limit to $M_\infty^4=M_0^4\cup{}M_0^4\cup\dots$ whose boundary is now entirely pseudoconvex.
		}} \\
	\end{tabular}
	\caption{
		Schematic depiction of the pseudoconvex manifold $M_\infty^4$.
	} \label{FigBuildingMInf}
\end{figure}

To remedy this, we create a sequence of manifolds-with-boundary $M_0^4\subset{}M_1^4\subset{}M_2^4\subset\dots$ that exhaust the entirety of a new pseudoconvex manifold $M_\infty^4$.
Set $M_0^4=M^4{}'$.
To create $M_1^4$, start with one copy of $M_0^4$, and glue another copy of $M_0^4$ onto the matching boundary or boundaries created by the pseudo-concave copy of $V_0^3$ (in case ({\it{ii}})) or the matching boundaries created by both $V_0^3$ and $V_1^3$ (in case ({\it{iii}})).
The manifold $M_1^4=M_0^4\cup{}M_0^4$ still has pseudo-concave boundary components: in case ({\it{ii}}) it has a pseudo-concave copy of $V_0^3$ and in case ({\it{iii}}) a pseudo-concave copy of both $V_0^3$ and $V_1^3$ (see the second image of Figure \ref{FigBuildingMInf}).
The gluing occurs by matching boundary components precisely as they had been matched in the original manifold and therefore the new manifold remains analytic, and continues to have, for example, an integrable complex structure.\\

For an induction process, starting with $M^4_{i-1}$ we attach a copy of $M_0^4$ to $M^4_{i-1}$ along the pseudo-concave copy of $V_0^3$ (or of both $V_0^3$ and $V_1^3$) to the matching pseudoconvex copy of $V_0^3$ (or of $V_0^3$ and $V_1^3$) of $M_0^4$.
Now we have K\"ahler manifolds-with-boundary $M_i^4=M_0^4\cup\dots\cup{}M_0^4$ (union of $i$ many copies of $M_0^4$), each $M_i^4$ created from $M_{i-1}^4$ by attaching one copy of $M_0^4$.

Picking any basepoint $p_0$, we now take a pointed limit in the Gromov-Hausdorff sense, and arrive at a K\"ahler manifold-with-boundary $M_\infty^4=M_0^4\cup{}M_0^4\cup\dots$.
This manifold has pseudoconvex boundary (the third image in Figure \ref{FigBuildingMInf}), but no pseudo-concave boundary components.
The new manifold $\partial{}M_\infty^4$ has exactly one pseudoconvex copy each of $V_0^3$ and $V_1^3$.

\underline{\it Proof that cases ({\it{ii}}) and ({\it{iii}}) cannot occur.}
By gluing to together copies of $M_0^4$, we have constructed the pseudoconvex manifold-with-boundary $M_\infty^4$ which has one pseudoconvex copy each of $V_0^3$ and $V_1^3$ on its boundary, and no pseudo-concave boundary components.
Both boundary components $V_0^3$, $V_1^3$ are compact, and one of the groups $\pi_1(V_0^3)$, $\pi_1(V_1^3)$ is finite.
This contradicts Lemma \ref{lemmaFundGrpsAndParaThm} and concludes the proof of Proposition \ref{PropEssSepPseudo}.

\vspace{0.1in}

\noindent{\it Proof of Theorem \ref{ThmPseudoconvexSubs}.}
The hypotheses are that $(M^4,J,g)$ has pseudoconvex boundary or no boundary, and $V^3$ is a boundary component or else a compact embedded submanifold that is pseudoconvex with respect to one of its two normals.
We also assume the group $\pi_1(V^3)$ is finite.

If $V^3$ is a boundary component, then Lemma \ref{lemmaFundGrpsAndParaThm} says it is the only boundary component and $\iota_*(\pi_1(V^3))=\pi_1(M^4)$.
If $V^3$ is separating, then $M^4\setminus{}V^3$ has a component $M_0^4{}'$ on which $V^3$ is pseudoconvex.
In this case, Lemma \ref{lemmaFundGrpsAndParaThm} provides the conclusion that $\partial{}M_0^4=V^3$ and $\iota_*\pi_1(V^3)=\pi_1(M_0^4)$.
The rest of the proof deals with the case that $V^3$ is non-separating.

Assuming $V^3$ is non-separating, we make a standard argument to show that $\pi_1(M^4)$ has a factor of $\mathbb{Z}$.
After removing $V^3$ we are left with a connected manifold $M_0^4=M^4\setminus{}V^3$, and the former submanifold $V^3$ becomes two distinct boundary components, one copy being pseudoconvex and the other being pseudo-concave.
For each $i\in\mathbb{Z}$ let $M_i^4$ be a copy of $M_0^4$.
Then gluing the pseudoconvex copy of $V^3$ in $M_i^4$ to the pseudo-concave copy of $V^3$ in $M_{i+1}^4$ and so on, we reach in the limit (a pointed Gromov-Hausdorff limit) the manifold
\begin{eqnarray}
	M_\infty^4\;=\;\dots\,\cup\,M_{-2}^4
	\,\cup\,M_{-1}^4\,\cup\,M_{0}^4\,\cup\,M_{1}^4\,\cup\,M_{2}^4
	\,\cup\,\dots \label{ContructionOfMInf}
\end{eqnarray}
Each $M_i^4$ is a copy of $M_0^4$, and each gluing surface maps unambiguously to $V^3\subset{}M^4$.
Therefore we retain a map $M_\infty^4\rightarrow{}M^4$.
Each point within any $M_i^4$---in particular any point on any of the gluing loci---has a small neighborhood which is evenly covered.
Therefore this is a covering map.
The map $d:M_\infty^4\rightarrow{}M_\infty^4$ generates the deck group of this covering.

We have constructed a cover of $M^4$ with transitive deck group $\left<d\right>\approx\mathbb{Z}$.
Covering space theory now says $\mathbb{Z}\approx\pi_1(M^4)/\pi_1(M^4_\infty)$, which is the same as the existence of a short exact sequence $0\rightarrow\pi_1(M^4_\infty)\rightarrow\pi_1(M^4)\rightarrow\mathbb{Z}\rightarrow0$.
But $\mathbb{Z}$ is projective so this sequence splits.
We conclude that $\pi_1(M^4)\approx\pi_1(M^4_\infty)\rtimes\mathbb{Z}$.

Next we verify that neither $M^4$ nor $M_\infty^4$ has any boundary components.
Consider the pseudoconvex manifold $\dots\cup{}M_{-2}^4\cup{}M_{-1}^4\cup{}M_0^4\subset{}M_\infty^4$.
This has a boundary component that is a pseudoconvex copy of $V^3$.
Because $V^3$ has $\pi_1(V^3)$ finite, Lemma \ref{lemmaFundGrpsAndParaThm} says this copy of $V^3$ is its \textit{only} boundary component.
If $M^4$ were itself to have a boundary component, this component would also be a component of the boundary of $M_0^4$.
Therefore $M_\infty^4$ would have additional boundary components.
We have just shown this is impossible, so we conclude $\partial{}M^4=\varnothing$.

Next we show $\pi_1(M^4_\infty)=\iota_*(\pi_1(V^3))$.
Consider the inclusion $\iota_*\pi_1(V^3)\subseteq\pi_1(M^4)\approx\pi_1(M_\infty^4)\rtimes\mathbb{Z}$.
Because $\pi_1(V^3)$ is a finite group, its inclusion must be an inclusion into the subgroup $\pi_1(M_\infty^4)\times\{0\}$---to see this, note that it is impossible for any element $a\in\pi_1(V^3)$ to map to any element of the form $(b,n)\in\mathbb{Z}\rtimes\pi_1(M_\infty^4)$ where $n\ne0$, for any such element would be infinite cyclic.
This produces an inclusion $\iota_*\pi_1(V^3)\subseteq\pi_1(M_\infty^4)$ (simply by leaving off the ``0'' factor).

For a proof by contradiction, assume the inclusion is strict: $\iota_*(\pi_1(V^3))\subset\pi_1(M_\infty^4)$.
Letting $p:\widetilde{M_\infty^4}\rightarrow{}M_\infty^4$ be the universal cover, the pre-image $p^{-1}(V^3)$ has more than one component.	
Referring to the construction of $M^4_\infty$, for any $k$ consider the submanifold
\begin{equation}
	M_k^4{}'
	\;\;\;\;=\;\;\;\;
	\dots\cup{}M_{k-2}^4\cup{}M_{k-1}^4\cup{}M_k^4
	\;\;\;\;\subset\;\;\;\;
	M_\infty^4.
\end{equation}
The submanifolds $\dots\subset{}M_{k-1}^4{}'\subset{}M_k^4{}'\subset\dots$ constitute an exhaustion of $M_\infty^4$.
Also note that $M_k^4{}'$ and $M_l^4{}'$ are diffeomorphic, by the deck action $d^{l-k}$.

Because $M_\infty^4$ has no boundary, as we saw above, the submanifold $M_k^4{}'$ has precisely one boundary component, which is a pseudoconvex copy of $V^3$.
The pre-image $\widetilde{M_k^4{}'}=p^{-1}(M_k^4{}')$ has more than one boundary component (because the pre-image of its boundary $p^{-1}(V^3)$ has more than one component).
However, conceivably $\widetilde{M_k^4{}'}$ also has more than one component, so we cannot immediately apply Lemma \ref{lemmaFundGrpsAndParaThm}.

So let $\gamma\subset{}M_\infty^4$ be a path that represents an element in $\pi_1(M^4_\infty)$ not in $\iota_*(\pi_1(V^3))$.
There exists some $k$ for which the loop $\gamma$ is contained within $M^4_k{}'$.
We may move the basepoint to the boundary copy of $V_0^3$, so without loss of generality we may assume $\gamma$ lies entirely within $M_k^4{}'$ and its two ends both terminate on the same point of the boundary copy of $V_0^3$.

Now lift both $M_k^4{}'$ and the path $\gamma$ along $p$.
Then the lifted path $\tilde{\gamma}$ is entirely inside $\widetilde{M_k^4{}'}$, and
because $\gamma$ represented an element of $\pi_1(M_\infty^4)$ not in $\iota_*(\pi_1(V^3))$, the path $\gamma$ no longer a loop but a path with distinct ends.
Because the basepoint in $M^4$ was originally on the boundary component $V_0^3$, the two ends of the lifted path now terminate on (the pseudo-concave sides of) two different boundary components.
Therefore at least one path-connected component of $\widetilde{M^4_k{}'}$ has two or more compact pseudoconvex copies of $V_0^3$ on its boundary.
These two boundary components are compact and have finite $\pi_1$, contradicting Lemma \ref{lemmaFundGrpsAndParaThm}.
Therefore $\iota_*(\pi_1(V^3))=\pi_1(M^4_\infty)$, as claimed.
We conclude that $\pi_1(M^4)\approx\iota_*(\pi_1(V^3))\rtimes\mathbb{Z}$.

Finally we prove that if $V^3$ is non-separating, then $M^4$ is both complete and non-compact.
We have already seen that the K\"ahler manifold $M^4$ has no boundary, so it is geodesically complete.
Because $\pi_1(M^4)=\iota_*(\pi_1(V^3))\rtimes\mathbb{Z}$ and $\iota_*(\pi_1(V^3))$ is finite, certainly $H_1(M^4;\mathbb{R})=\mathbb{R}$, so the first betti number of $M^4$ is one.
Thus it is impossible that $M^4$ is compact, because the first betti number is even on all compact K\"ahler manifolds.
This concludes the proof of Theorem \ref{ThmPseudoconvexSubs}.

\section{Boundaries, holomorphic submanifolds, and manifold ends}

In this section we prove the four corollaries, but first we require some information on neighborhoods of embedded holomorphic $\mathbb{P}^1$ submanifolds of non-negative self-intersection.
From the differentiable perspective such neighborhoods are completely canonical: an embedded $\mathbb{P}^1$ has a neighborhood diffeomorphic to its normal bundle.
The normal bundle is a complex line bundle over $\mathbb{P}^1$, and differentiably these are fully characterized by the self-intersection number of the zero-section.

But to apply our theorems we require more: neighborhood boundaries must have a sign on their Levi forms.
In the case of zero self-intersection, an embedded $\mathbb{P}^1$ has a neighborhood that is a holomorphic product of $\mathbb{P}^1$ with a disk in $\mathbb{C}$ \cite{Sav82}.
Certainly the boundary of such a neighborhood is Levi-flat and has the topology of $\mathbb{S}^2\times\mathbb{S}^1$.
In the case of \textit{positive} self-intersection the situation is known to be much more complicated: neighborhoods of $\mathbb{P}^1$ of positive self-intersection fall within an infinite-dimensional moduli space of complex structures \cite{Mish93}.
However all we require is pseudo-concavity, and pseudo-concave neighborhoods are indeed known to exist.
The following result, an important thoerem in its own right, we record for our purposes as a lemma.
\begin{lemma}[cf. \cite{AG62} \cite{Rossi64} \cite{Suz75} \cite{Sav82}] \label{LemmaPosCycle}
	Let $(M^4,J)$ be a complex manifold and assume $\mathbb{P}^1$ is an embedded holomorphic $\mathbb{P}^1$ of positive self-intersection.
	Then a neighborhood $\Omega$ of this $\mathbb{P}^1$ exists, whose boundary has negative Levi-form with respect to the outward pointing normal.
	Further, $\Omega$ deformation-retracts onto $\mathbb{P}^1$, $\partial\Omega$ is diffeomorphic to a lens space, and we can make this neighborhood as small as desired: given any open set $\Omega'$ containing $\mathbb{P}^1$, we can assume $\overline{\Omega}\subset\Omega'$.
	
	Likewise, if such a $\mathbb{P}^1$ has zero self-intersection, there exist a small neighborhood $\Omega$ of $\mathbb{P}^1$ with Levi-flat boundary and boundary topology $\mathbb{S}^2\times\mathbb{S}^1$.
\end{lemma}

{\it Proof of Corollary \ref{CorLeviFlatAndCycle}.}
The assumption is that $(M^4,J,g)$ has a compact Levi-flat surface $V^3$, as well as a $\mathbb{P}^1$ of positive self intersection.

If $\mathbb{P}^1$ and $V^3$ do not intersect, then using the metric we can build a small tubular neighborhood around $\mathbb{P}^1$ that does not intersect $N^3$.
From Lemma \ref{LemmaPosCycle}, inside this tubular neighborhood is a domain $\Omega$ containing $\mathbb{P}^1$ that is pseudo-concave and compact, and $\partial\Omega$ has finite fundamental group.
Thus $M^4\setminus\Omega$ has a pseudo\textit{convex} boundary component.
The boundary of $M^4\setminus\Omega$ is compact, pseudoconvex, and has finite fundamental group, so Theorem \ref{ThmPseudoconvexSubs} says $\pi_1(M^4)$ is finite.

Therefore $(M^4\setminus\Omega)\setminus{}V^3$ is disconnected (by an easy Mayer-Vietoris argument).
Let $M^4{}'$ be the component of $(M^4\setminus\Omega)\setminus{}V^3$ whose boundary $\partial{}M^4{}'$ is the union of the compact pseudoconvex submanifold $\partial\Omega$, which has finite fundamental group, and a copy of $V^3$, which is compact and Levi-flat.
This contradicts Thoerem \ref{ThmPseudoconvexSubs} and establishes the Corollary.

{\it Proof of Corollary \ref{CorCycsPseudo}.}
This corollary studies two possibilities: that a holomorphic embedded $\mathbb{P}^1$ of positive self-intersection exists, and that a holomorphic embedded $\mathbb{P}^1$ of zero self-intersection exists.

First assume $(M^4,J,g)$ has a $\mathbb{P}^1$ of positive self-intersection.
By Lemma \ref{LemmaPosCycle} we can place a sufficiently small neighborhood $\Omega$ around $\mathbb{P}^1$ so that $\Omega$ does not intersect $\partial{}M^4$ and so $M^4\setminus\Omega$ is still pseudoconvex.
The boundary component of $M^4\setminus\Omega$ corresponding to $\partial\Omega$ is a sphere or lens space, so has finite fundamental group.
Because there is one compact pseudoconvex boundary component with finite fundamental group, by Theorem \ref{ThmPseudoconvexSubs} this is the \textit{only} component of $\partial(M^4\setminus\Omega)$.
Therefore $\partial{}M^4=\varnothing$.

Still assuming $\mathbb{P}^1$ has positive self-intersection, assume $\pi_1(M^4)$ is not trivial.
Then the universal cover $\widetilde{M}{}^4$ of $M^4$ is non-trivial, so there exists two (or more) embedded holomorphic $\mathbb{P}^1$ submanifolds of positive self-intersection.
By Lemma \ref{LemmaPosCycle}, around these we can find non-intersecting pseudo-concave neighborhoods $\Omega_1$ and $\Omega_2$ so that $\partial\Omega_1$, $\partial\Omega_2$ both have finite fundamental groups.
Then $\widetilde{M}{}^4\setminus(\Omega_1\cup\Omega_2)$ has pseudoconvex boundary and two of its compact components have finite fundamental group.
This violates Theorem \ref{ThmPseudoconvexSubs}.

To verify the second assertion, assume some $\mathbb{P}^1$ has zero self-intersection.
For a contradiction assume there is a compact boundary component $V^3\subseteq\partial{}M^4$ with finite $\pi_1(V^3)$.
By Lemma \ref{LemmaPosCycle}, $\mathbb{P}^1$ has a small neighborhood $\Omega$ that does not intersect $V^3$, and is Levi-flat and compact.
However because $\pi_1(V^3)$ is finite, the manifold $M^4\setminus\Omega$ violates Theorem \ref{ThmPseudoconvexSubs}.
This establishes the corollary.

{\it Proof of Corollary \ref{CorZeroAndPos}.}
The assumption is that $(M^4,J,g)$ is K\"ahler manifold, without boundary (not necessarily compact), that $N=\mathbb{P}^1$ has non-negative self-intersection, and that $N'=\mathbb{P}^1$ has positive self-intersection..

For a proof by contradiction assume $N\cap{}N'=\varnothing$.
By Lemma \ref{LemmaPosCycle}, $N$ and $N'$ have small, non-intersecting neighborhoods $\Omega$ and $\Omega'$; the boundary of $\Omega$ is Levi-flat and compact and the boundary of $\Omega'$ is compact, pseudo-concave, and has finite fundmantal group.
Therefore $M^4\setminus(\Omega\cup\Omega')$ has two compact pseudoconvex boundary components, one of which has finite fundamental group.
This violates Theorem \ref{ThmPseudoconvexSubs}, and establishes the corollary.

{\it Proof of Corollary \ref{CorEnds}.}
The assumption is that $(M^4,J,g)$ has $k$ many ALE ends and $l$ many ALF ends; apriori $k$ or $l$ might be infinite.

Any ALE or ALF end has a geometrically convex separating surface (by ``geometrically convex'' we mean non-negative second fundamental form), and because these ends are geometrically convex they are pseudoconvex (this is because, in the K\"ahler setting, the Levi form is the $J$-averaging of the second fundamental form, restricted to the distribution perpendiculat to the normal $\hat{n}$ and $J\hat{n}$).
Label these separating submanifolds $V_1^3,V_2^3,\dots$; there might be finitely or infinitely many of these.
Each $V_i^3$ separates the manifold into two pieces: a manifold end (with a single boundary component, which is pseudo-concave), and a manifold with a new pseudoconvex boundary component.
The boundary components $V_i^3$ are each quotients of spheres, so have finite fundmanental group.
If there are two such ends, then a component of $M^4\setminus(V_1^3\cup{}V_2^3)$ has two pseudoconvex compact boundary components with finite fundamental group, in violation of Theorem \ref{ThmPseudoconvexSubs}.
Therefore at most one such end exists, so $k+l$ is zero or one.

Next we suppose $(M^4,J,g)$ has an embedded holomorphic $\mathbb{P}^1$ submanifold with non-negative self-intersection.
By Lemma \ref{LemmaPosCycle} there is a neighborhood $\Omega$ so that $M^4\setminus\Omega$ has a compact, pseudoconvex boundary.
But then if $M^4$ has an ALE or ALF end, we can separate that end, and have a compact pseudoconvex boundary remaining, that is topologically either $\mathbb{S}^3$ or a Lens space so has finite fundmantal group.
There are now two compact pseudoconvex boundary components and one has finite fundamental group, contradicting Theorem \ref{ThmPseudoconvexSubs}.

Lastly if $k+l=1$, $\pi_1(M^4)$ is finite by Lemma \ref{lemmaFundGrpsAndParaThm}.
This establishes the corollary.

\section{Examples} \label{SecExamples}

\subsection{A K\"ahler surface with non-intersecting Levi-flat submanifolds}

A K\"ahler $(M^4,J,g)$ cannot have Levi-flat submanifolds $V_0^3$, $V_1^3$ with $V_0^3\cap{}V_0^3=\varnothing$ when either submanifold has finite fundamental group (by Theorem \ref{ThmLeviFlats}).
We give an example showing that if both submanifold have infinite fundmantal group, then two such submanifolds can exist.

Let $M^4$ be the complex surface $\mathbb{P}^1\times\mathbb{P}^1$.
Giving this the usual product metric, this is a K\"ahler manifold.
Pick two points $p,q\in\mathbb{P}^1$, and consider the complex submanifolds $\{p\}\times\mathbb{P}^1$ and $\{q\}\times\mathbb{P}^1$.
Expand these into small neighborhoods
\begin{eqnarray}
	\Omega_\epsilon\;=\;D_\epsilon(p)\times\mathbb{P}^1, \quad
	\Omega_\epsilon'\;=\;D_\epsilon(q)\times\mathbb{P}^1
\end{eqnarray}
where the disks are defined as follows: Let $f:\mathbb{P}^1\setminus\{p,q\}\rightarrow\mathbb{C}$ be the harmonic function with a simple pole of $-\infty$ at $p$ and a simple pole of $+\infty$ at $q$.
Then define $D_\epsilon(p)=\{f<\log(\epsilon)\}$ and $D_\epsilon(q)=\{f>-\log(\epsilon)\}$.

Since $f$, interpreted as an extended real-valued function $f:\mathbb{P}^1\times\mathbb{P}^1\rightarrow\mathbb{R}\cup\{-\infty,\infty\}$ is harmonic on one factor and constant on the other, it is pluriharmonic.
Consequently the boundaries $\partial\Omega_\epsilon=\{f=\log\epsilon\}$ and $\partial\Omega_\epsilon'=\{f=-\log\epsilon\}$ are both Levi-flat.
For all sufficiently small $\epsilon$ these neighborhoods are non-intersecting.
They are both topologically equivalent to $\mathbb{S}^2\times\mathbb{S}^1$, so have infinite fundamental group.

\subsection{Levi-flat surfaces in $\mathbb{P}^2\#\overline{\mathbb{P}}{}^2$} \label{SubsecRationalExample}

In addition to the example in $\mathbb{P}^1\times\mathbb{P}^1$, we give examples of Levi-flats within $M^4=\mathbb{P}^2\#\overline{\mathbb{P}}{}^2$.
This is a rational surface, and has holomorphic $\mathbb{P}^1$ submanifolds of self-intersection $+1$.

In $\mathbb{C}^2=\{(z_1,z_2)\}$ an example of a Levi-flat hypersurface is $\{Im(z_1)=\alpha\}$ where $\alpha$ is any fixed real number.
Embedding this into $\mathbb{P}^2$ by $(z_1,z_2)\mapsto[1;z_1;z_2]$ and then taking the closure produces a closed Levi-flat hypersurface $V^3_\alpha$, except that $V_\alpha^3$ has a singular point which lies on the ``sphere at infinity'' $\mathbb{P}^1=[0;z_1;z_2]$.
We examine how $V_\alpha^3$ interacts with the curves of positive self-intersection.
One easily checks that the $\mathbb{P}^1$ at infinity entirely lies within $V^3_\alpha$, which also contains the singular point.
The other rational curves of positive self-intersection are (closures in $\mathbb{P}^2$ of) the 1-dimensional subspaces within $[1;z_1;z_2]$.
Each of these intersects $V_\alpha^3$ either in a linear subspace and
lies completely within $V_\alpha^3$ (these are the various $\mathbb{P}^1$ given by $z_1=const$ where $Im(const)=\alpha$), or in a circle (these are the various $\mathbb{P}^1$ given by the linear subspaces of $\mathbb{C}^2$ \textit{not} of the form $z_1=const$), or only at its singular point (these are the various $\mathbb{P}^1$ given by any subspace of the form $z_1=const$ where $Im(const)\ne\alpha$).
The space $V_\alpha^3$ is foliated by the complex lines $z_1=const$, $Im(const)=\alpha$, and each of these intersects the singular point.

Next we look at the singular point more closely.
As usual there are three natural charts for $\mathbb{P}^2$.
We defined $V_\alpha^3$ in the chart $(z_1,z_2)\mapsto[1;z_1;z_2]$ with $V_\alpha^3$ being $\{[1,\alpha'+\sqrt{-1}\alpha;z_2]\}$, $\alpha'\in\mathbb{R}$, $z_2\in\mathbb{C}$.
In the chart $[z_0;1;z_2]$ the surface $V_\alpha^3$ appears as $Im(1/z_0)=\alpha$, which remains non-singular and in this chart is a copy of $\mathbb{S}^1\times\mathbb{R}^2$ (except for the case $\alpha=0$).
Finally in the chart $[z_0;z_1;1]$ the hypersurface $V_\alpha^3$ appears as $Im(z_1/z_0)=\alpha$, which is a complex cone at $(0,0)$.
Since $z_1/z_0$ is invariant under simultaneous multiplication of both variables, we see a foliation of this cone by complex lines.
The lines themselves are parameterized by $\mathbb{S}^1$: this is because with $Im(z_1/z_0)$ fixed, the complex line on which $(z_0,z_1)$ lies is determined by $Re(z_1/z_0)\in\mathbb{R}\cup\{\infty\}\approx\mathbb{S}^1$.

From the singular surface $V_\alpha^3$ we create a \textit{smooth} Levi-flat hypersurface which is diffeomorphic to $\mathbb{S}^2\times\mathbb{S}^1$.
Blow up $\mathbb{P}^2$ at the singular point of $V_\alpha^3$.
In particular, this desingularizes $V_\alpha^3$, as its cone point was removed, and replacing it with a copy of $\mathbb{S}^1$ (due to the fact that the complex lines of $V_\alpha^3$ passing through the cone point were paramterized by $\mathbb{S}^1$).
One can check on charts that the new hypersurface, which is now within $\mathbb{P}^2\#\overline{\mathbb{P}}^2$, is smooth.
Let $V_\alpha^3{}^*$ be this new, non-singular surface.

As we have seen, the intersection of the desingularized surface $V_\alpha^3{}^*$ with the exceptional divisor is a circle.
On $V^3_\alpha{}^*$ itself, the desingularization can be seen by removing a neighborhood of the singular point---the boundary of this neighborhood is a torus---and then gluing in a solid torus whose central circle is the intersection locus with the exceptional divisor.
The boundary torus is glued by a torus automorphism---the automorphism exchanges a diagonal circle with a contractible meridian and vice-versa---and so the desingularized surface has infinite fundamental group.
This surface continues to intersect every rational curve of positive self-intersection.
It is now foliated by curves of zero self-intersection.

The original singular hypersurfaces $V_\alpha^3\subset\mathbb{P}^2$, which are simply connected, are themselves of interest.
These all intersect every $\mathbb{P}^1$ of positive self-intersection, so each $V_\alpha^3$ obeys Corollary \ref{CorLeviFlatAndCycle}, even though they are not $C^2$.
The singular hypersurfaces $V_\alpha^3$ all intersect one another (along the $\mathbb{P}^1$ at infinity), so the conclusion of Theorem \ref{ThmLeviFlats} also holds for the surfaces $V_\alpha^3$.

In the future, one would like to extend some version of Theorem \ref{ThmLeviFlats} to the case of potentially singular Levi-flat surfaces, as singular Levi-flats seem to be plentiful.
A number of technical hurdles stand in the way, such as questions of internal topology, possible nodes or self-intersections, how the Levi-flat condition interacts with or persists across the singularity, and points of contact between submanifolds possibly being at singular points.

See \cite{N99} for additional examples of Levi-flat submanifolds in ruled surfaces.

\subsection{A scalar-flat 2-ended K\"ahler instanton as a Taub-NUT conformal transform} \label{SubSecTwoEnded}

Our final example shows that a K\"ahler 4-manifold, even a scalar-flat K\"ahler 4-manifold, can have two ends even if one end is ALE.

The classic Euclidean Taub-NUT metric of Hawking \cite{Haw77} is a startling example in a number of ways, one of these ways being that although K\"ahler itself, it is conformal to two very different K\"ahler metrics.
Normally it is impossible that a K\"ahler metric in dimension 4 be conformal to another K\"ahler metric, but here it is possible due to the unusually large number of complex structures the Taub-NUT metric is compatible with.
The classic metric is
\begin{equation}
	g\;=\;\frac14\frac{r+m}{r-m}dr^2
	+4m^2\frac{r-m}{r+m}(\sigma^1)^2
	+(r^2-m^2)\big((\sigma^2)^2+(\sigma^3)^2\big)
\end{equation}
on $r\in[m,\infty)$ where $m$ is a constant, and $\sigma^1$, $\sigma^2$, $\sigma^3$ are the usual left-invariant 1-forms on $\mathbb{S}^3$ normalized in the usual way so $d\sigma^i=-\epsilon^i{}_{jk}\sigma^j\wedge\sigma^k$.
It is well known that this metric is hyperK\"ahler.
It is half-conformally flat, but in the usual orientation $\{dr,\sigma^1,\sigma^2,\sigma^3\}$ the ``wrong'' half: $s=0$ and $Ric=0$ but
\begin{equation}
	W^+\;=\;\frac{8m}{(r+m)^3}
	\left(\frac32\omega\otimes\omega
	-Id_{\bigwedge^-} \right),
	\quad\quad W^-\;=\;0 \label{EqnTaubNUTDef}
\end{equation}
(in those days physicists were seeking ``self-dual'' metrics).
The hyperK\"ahler structure is carried by three sections $\omega_1,\omega_2,\omega_3\in\bigwedge^-$ which dualize to compatible complex structures $I_1$, $I_2$, $I_3$ that obey the quaternionic relations (and give the orientation opposite to $\{dr,\sigma^1,\sigma^2,\sigma^3\}$).

However the metric (\ref{EqnTaubNUTDef}) is compatible with two additional complex structures, not related to $I_1$, $I_2$, or $I_3$, which we label $J^+$ and $J^-$.
Letting $e_1,e_2,e_3$ be the left-invariant fields on $\mathbb{S}^3$ dual to the forms $\sigma^i$, and given any function $f=f(r)$, let $J_f$ be the almost complex structure
\begin{equation}
	J_f\;=\;
	-f\frac{\partial}{\partial{}r}\otimes\sigma^1
	+\frac{1}{f}e_1\otimes{}dr
	-e_2\otimes\sigma^3
	+e_3\otimes\sigma^2.
\end{equation}
In fact this is always integrable.
To see this, note that
\begin{equation}
	\bigwedge{}^{0,1}=span_{\mathbb{C}}
	\left\{\frac{1}{f}dr-\sqrt{-1}\sigma^1,\;\sigma^2-\sqrt{-1}\sigma^3
	\right\},
\end{equation}
and then on the basis we compute
\begin{equation}
	\small
	\begin{aligned}
		&d\left(\frac{1}{f}dr-\sqrt{-1}\sigma^1\right)
		=2\sqrt{-1}\sigma^2\wedge\sigma^3
		=-2\sigma^2\wedge\left(\sigma^2-\sqrt{-1}\sigma^3\right) \\
		&d\left(\sigma^2-\sqrt{-1}\sigma^3\right)
		\;=\;2\sigma^1\wedge\sigma^3
		+2\sqrt{-1}\sigma^1\wedge\sigma^2
		\;=\;2\sqrt{-1}\sigma^1\wedge\left(\sigma^2-\sqrt{-1}\sigma^3\right).
	\end{aligned}
\end{equation}
Therefore $d\left(\bigwedge^{0,1}\right)\subset\bigwedge^1\wedge\bigwedge^{0,1}$ so $J_f$ is integrable.
Then two complex structures compatible with $g$ are
\begin{equation}
	\begin{aligned}
		&J^+\;=\;
		-4m\frac{r-m}{r+m}\frac{\partial}{\partial{}r}\otimes\sigma^1
		+\frac{1}{4m}\frac{r+m}{r-m}e_1\otimes{}dr
		-e_2\otimes\sigma^3
		+e_3\otimes\sigma^2 \\
		&J^-\;=\;
		\;\;\,4m\frac{r-m}{r+m}\frac{\partial}{\partial{}r}\otimes\sigma^1
		-\frac{1}{4m}\frac{r+m}{r-m}e_1\otimes{}dr
		-e_2\otimes\sigma^3
		+e_3\otimes\sigma^2.
	\end{aligned}
\end{equation}
Their corresponding $(1,1)$ forms are
\begin{equation}
	\begin{aligned}
		&g\left(J^+\cdot,\,\cdot\right)
		\;=\;\;\;\,m\,dr\wedge\sigma^1+(r^2-m^2)\sigma^2\wedge\sigma^3\;\in\;\bigwedge{}^+, \\
		&g\left(J^-\cdot,\,\cdot\right)
		\;=\;-m\,dr\wedge\sigma^1+(r^2-m^2)\sigma^2\wedge\sigma^3\;\in\;\bigwedge{}^-.
	\end{aligned}
\end{equation}
Using the facts that $d\sigma^1=-2\sigma^2\wedge\sigma^3$ and $d(\sigma^2\wedge\sigma^3)=0$ we easily compute
\begin{equation}
	d\left(g\left(J^\pm\cdot,\,\cdot\right)\right)
	\;=\;d\big(\log(r\mp{}m)^2\big)\wedge{}g\left(J^\pm\cdot,\,\cdot\right),
\end{equation}
so neither form is K\"ahler, but both are conformally K\"ahler.
Making the conformal change with factors $(r\mp{}m)^2$ gives the two new metrics
\begin{equation}
	g^{\pm}\;=\;(r\mp{}m)^{-2}g.
\end{equation}
The corresponding forms $\omega^\pm=g^\pm(J^\pm\cdot,\cdot)$ are now both closed, so we have two new K\"ahler metrics $g^+$ and $g^-$.
Both are conformally related to $g$ so they retain half-conformal flatness.
In particular they are Bach-flat, so by Proposition 4 part (v) of \cite{Derd}, they are both extremal K\"ahler metrics.
The pair of K\"ahler structures $\{(g^+,J^+,\omega^+),(g^-,J^-,\omega^-)\}$ on $M^4$ is an \textit{ambiK\"ahler} structure; see \cite{ACG16}.

The metric $g^+$ has K\"ahler form $\omega^+\in\bigwedge{}^+$ and by Derdzinski's Theorem, Proposition 2 of \cite{Derd}, we have that $W^+=\frac{s}{12}\left(\frac32\omega^+\otimes\omega^+-Id_{\bigwedge^+}\right)$.
Therefore (\ref{EqnTaubNUTDef}) and the conformal invariance of $W^+$ shows that $s_{g^+}=\frac{96m}{r+m}$ and we see this metric is strictly extremal (meaning it is extremal but its scalar curvature is not constant).
Clearly $g^+$ is one-ended, as the conformal factor $(r+m)^{-2}$ is non-singular.
One easily checks $g^+$ is complete.
However, unlike the Taub-NUT, it is not ALF.
Asymptotically it is cusp-like: $g^+\approx\frac14r^{-2}(dr^2+(4m)^2(\sigma^1)^2)+(\sigma^2)^2+(\sigma^3)^2$, which as $r\rightarrow\infty$ approaches the metric product of a pseudo-sphere times a sphere (locally, not globally).

The $g^-$ metric is more interesting to us, as it is two-ended.
It has K\"ahler form $\omega^-\in\bigwedge^-$.
Since $W^-=0$, Derdzinski's theorem states $g^-$ is scalar-flat.
Alternatively, one could use the conformal change formulas to compute $s=0$.
In addition to $s=0$ and $W^-=0$, the conformal change formulas show
\begin{equation}
	\begin{aligned}
		&Ric_{g^-}=4\left(\frac{r-m}{r+m}\right)^2\left(
		-(\eta^0)^2-(\eta^1)^2+(\eta^2)^2+(\eta^3)^2
		\right) \\
		&W^+_{g^-}=8\frac{(r-m)^2}{(r+m)^3}\left(\frac{3}{2}|\omega|_{g^-}^{-2}\omega\otimes\omega
		-Id_{\bigwedge^+}\right)
	\end{aligned}
\end{equation}
where $\eta^0=|dr|^{-1}_{g^-}dr$ and $\eta^i=|\sigma^i|^{-1}_{g^-}\sigma^i$ are unit 1-forms.
In particular, all curvature components are uniformly bounded.
Notice that the Ricci curvature does not decay to 0 along the cusp-like end (as $r\rightarrow\infty$), but does decay to zero along the AE end (as $r\searrow{}m$).

However $g^-$ is two-ended, as the conformal factor $(r-m)^{-2}$ is singular at the ``nut'' at $r=m$.
In dimension 4, it is well-known that an inverse-quadratic conformal factor removes a point and replaces it with an asymptotically Euclidean end; see Theorem 6.5 of \cite{LP87} for example.
In addition to this new AE end that replaced the nut, there remains the end at $r=\infty$.
This end, in the $g^-$ metric, is asymptotically identical to the $g^+$ end: cusp-like, in the sense that it is (locally) very close to the metric product of the sphere with the psuedosphere.

Thus the metric $g^-$, conformally related to the Taub-NUT metric, is a scalar-flat K\"ahler 2-ended metric, with one asymptotically Euclidean end and one cusp-like end.



\end{document}